\documentclass[10pt]{amsart}

\usepackage{enteteanglais}

\title{On the entropy of Hilbert Geometries of Low Regularities}

\author{Jan Cristina and Louis Merlin}

\thanks{}
\begin{document}


\begin{abstract}
We compare the regularity of the boundary of a convex set with the value of its Finslerian volume entropy. The main result states that the volume entropy of a two-dimensional domain whose associated curvature measure is Ahlfors $\alpha$-regular is $\frac{2\alpha}{\alpha +1}$.
\end{abstract}

\maketitle

\section{Introduction}

\subsection{Hilbert geometries.}

To define a Hilbert geometry, we need a convex compact subset $\Omega$ of $\mathbb{R}^n$ (or a strictly convex set in $\mathbb{R}\mathbb{P}^n$). Then we construct a distance in the interior of the convex set using the cross-ratio.

Precisely, take two points $p$ and $q$ in Int$(\Omega)$. The compactness and convexity show that there exist two uniquely determined points $a$ and $b$ on $\partial \Omega$ such that $a$, $p$, $q$, $b$ are aligned in this order. We set
\[d_\Omega(p,q)=\frac{1}{2}\abs{\log\frac{\abs{q-a}\abs{p-b}}{\abs{p-a}\abs{q-b}} }. \]

\begin{figure}[ht]
\begin{center}
\includegraphics[scale=.3]{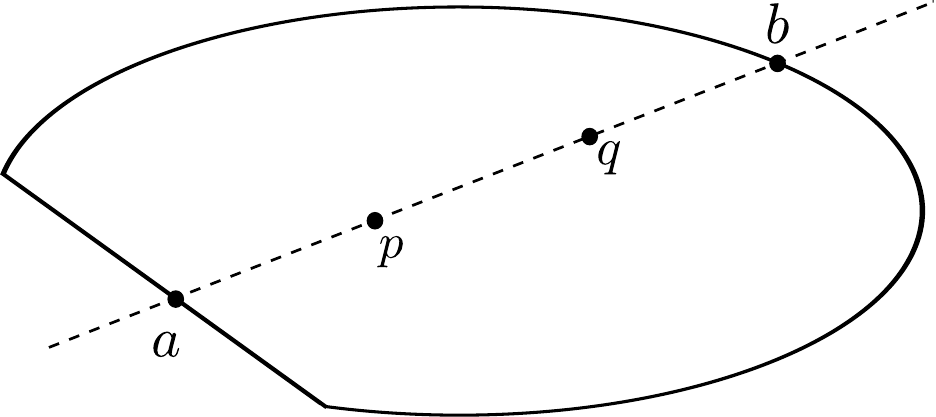}
\caption{Hilbert metric}
\end{center}
\end{figure}

\begin{figure}[ht]
\begin{center}
\includegraphics[scale=.3]{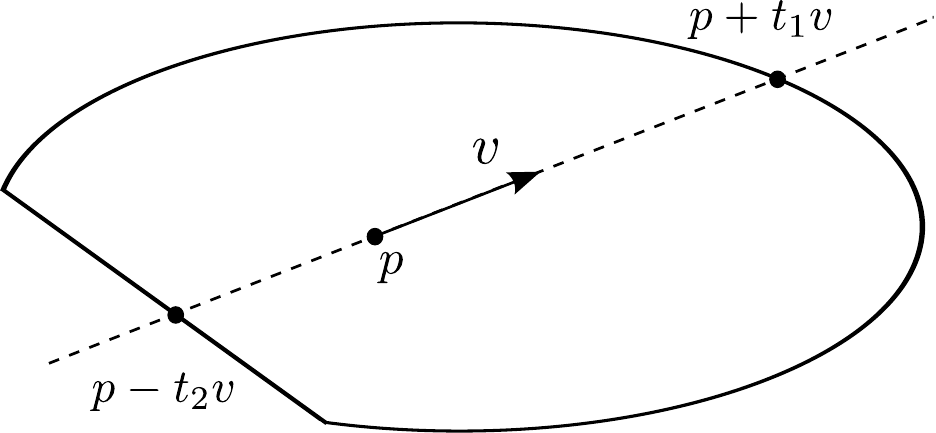}
\caption{Finslerian structure}
\label{fig:finsler}
\end{center}
\end{figure}

When $\Omega$ is an ellipsoid, we construct the Klein model for hyperbolic geometry. For any other case, the distance is not even Riemannian (\cite{Kay}). However it is Finslerian, infinitesimally generated by the field of norms $\left(\norm{\cdot}_x\right)_{x\in \Omega}$ given by
\[\norm{v}_x = \frac{1}{2}\left(\frac{1}{t_1}+\frac{1}{t_2} \right) \]
where $t_1$ and $t_2$ are positive real numbers such that $x+t_1v$ and $x-t_2v$ meet the boundary of $\Omega$ (figure \ref{fig:finsler}).

The couple $(\Omega,d_\Omega)$ is called a Hilbert geometry. Such spaces enjoy the following remarkable properties \cite{handbookhilbert}:
\begin{itemize}
\item The metric spaces $(\Omega,d_\Omega)$ are complete.
\item (Affine or projective) straight lines are geodesics.
\item The group of projective transformations leaving $\Omega$ invariant acts by isometries on $(\Omega,d_\Omega)$.
\end{itemize}

\begin{rmk}\label{convention}
Throughout this text, we will use the notation $\abs{\cdot}$ for refering to Euclidean lenghts, norms or volumes and $\norm{\cdot}$ will be used for the Finslerian quantities. 
\end{rmk}

\subsection{Entropy.}

In this context of Hilbert geometries, the main goal of this article is to study an invariant, the \textit{volume growth entropy}.
The volume entropy of a metric measured space $(X,d,\mu)$ is the exponential asymptotic growth rate of volume of balls when the radius goes to infinity. Precisely, it is the real number defined as the limit (whenever it exists) of 
\[\frac{\log (\mu (B(x,R)))}{R} \]
(it is independent of the choice of a basepoint $x$).
It is known to be a powerful invariant. When $(X,d)$ is a Riemannian manifold and $\mu$ is the Riemannian volume, it has been used several times to capture a lot of informations about the ambient geometry (see for instance \cite{bcg2}).

Returning to the context of a Hilbert geometry $(\Omega,d_\Omega)$ given by a convex $\Omega\subset \mathbb{R}^n$, one of the common issue is that the volume is not canonically defined. However we can isolate axioms for a notion of \textit{appropriate volume} (\cite{hilbertvolume} for which the entropy does not depend on the choice of an appropriate volume. In this paper, we compute volumes with the so-called Hausdorff ($n$-)measure. Let us define it.

We consider the function $\sigma$ on $\Omega$ given by
\[\sigma(x)=\frac{\omega_n}{\mathcal{L}(B_x)}\]
where $\omega_n$ is the measure of the unit Euclidean ball of $\mathbb{R}^n$, $\mathcal{L}$ is the Lebesgue measure and $B_x$ is the Finslerian unit ball. Finally the Hausdorff measure $\mu$ is the measure (absolutely continuous with respect to the Lebesgue measure), the density of which is given by $\sigma$, i.e:
\[\mu(A)=\int_A\sigma(x)d\mathcal{L}(x)\]
for any Borel set $A$. The density $\sigma$ is called the \textit{Busemann function}. Tha fact that $\mu$ is indeed the $n$-dimensional Hausdorff measure follows from \cite{busemannvolume}.

For a general Hilbert geometry, the quantity $\displaystyle \frac{\log (\mu (B(x,R)))}{R}$ does not converge in general when $R$ goes to infinity (\cite{entropyapprox} corollary 4), so we usually consider lower and upper entropies:
\[ \underline{h}(\Omega) = \liminf_{R\to\infty} \frac{\log (\mu (B(x,R)))}{R} \quad \mbox{ and }\quad\overline{h}(\Omega)=\limsup_{R\to\infty} \frac{\log (\mu (B(x,R)))}{R}         .\]
For instance, when $\Omega$ is an ellipsoid in $\mathbb{R}^n$,
\[ \underline{h}(\Omega)=\overline{h}(\Omega)= n-1.\]
The notation $h$ implicitly means that $\displaystyle \frac{\log (\mu (B(x,R)))}{R}$ converges.

The volume entropy of Hilbert geometries has been studied by various authors. For instance it has been proved by N. Tholozan in \cite{tholozanhilbert} that the entropy of Hilbert geometries never exceed the hyperbolic entropy. This result has been known in dimension 2 and 3 since the work of C. Vernicos in \cite{entropyapprox} and in the case of divisible convex sets thanks to a result of M. Crampon in \cite{cramponentropydivisible}. It is known that the entropy can take any value in dimension 2 (\cite{entropyapprox} corollary 4) and it has been precisely computed in many cases: for instance the entropy vanishes for convex polytopes \cite{entropypolytopes} and it is extremal as soon as the convex set is sufficiently regular. Actually the latter statement is the starting point of this paper, let us make it precise:

\begin{thm}[First main theorem in \cite{bbvc11}]\label{vernicos}
Suppose the boundary of the convex set $\Omega$ is a hypersurface of $\mathbb{R}^n$ of regularity $\mathcal{C}^{1,1}$. Then the entropy exists and
\[h(\Omega)=n-1. \]
\end{thm}

\subsection{Question.}

The question we would like to address in this paper is suggested by the previous result and is the following

\begin{quest}
Can we find a relation between the regularity of the boundary of a convex domain and the value of its Finslerian volume entropy ?
\end{quest}

In particular, is there a lower bound for the entropy of a domain given the regularity of its boundary? And are there lower regularities for the boundary than $C^{1,1}$ which guarantee maximal entropy?

\subsection{Results.}

We propose two types of answer for the previous question.
We first describe a bijective relation between entropy and regularity in dimension 2. Following theorem \ref{vernicos}, a natural (but slightly naive) approach would be to try to understand a $\mathcal{C}^{1,\alpha}$ regular convex set. As a global invariant, the volume entropy only sees the part of the convex set with maximal volume growth. In particular, as soon as there is a part of the boundary of $\Omega$ which has regularity $\mathcal{C}^2$ with positive Gauss curvature, then $h(\Omega)=n-1$. This may happen for a $\mathcal{C}^{1,\alpha}$ regular convex set and this implies that in particular there is no hope of characterizing the entropy of a domain with its regularity $\mathcal{C}^{1,\alpha}$. Hence, in order to precisely compare the regularity and the entropy, we need the regularity to realize the following constraint: "Whenever the Gauss curvature is positive (in a possibly weak sense), then the regularity must be strictly $\mathcal{C}^{1,\alpha}$ but not more." This is roughly the definition of Ahlfors $\alpha$-regularity. We postpone until section \ref{sec:prel2} the precise definition of Ahlfors regularity. The main result of this paper is then the following.

\begin{thm}[First main theorem]\label{thm:ahlfors}
Let $\Omega$ be a convex and relatively compact domain of $\mathbb{R}^2$ which is Ahlfors $\alpha$-regular. Then
\[h(\Omega)=\frac{2\alpha}{\alpha + 1}.\]
\end{thm}

In a second and independent part of this paper, we show a stronger version of theorem \ref{vernicos}, weakening the assumption of $C^{1,1}$-regularity. The space of $C^{1,1}$ maps is isomorphic to the Sobolev space $W^{2,\infty}$.  We show that for some finite values of $p$, domains whose boundaries are $W^{2,p}$-regular must have volume entropy equal to $n-1$. \emph{i.e.} we show that we can decrease the Sobolev regularity and still be sure to get a convex set of maximal volume entropy.

\begin{thm}[Second main theorem]\label{thm:sobolev}
  Let $\Omega$ be a convex relatively compact subset of $\mathbb{R}^n$. Assume that the boundary is parametrized by a homeomorphic map $\varphi : S^{n-1}\rightarrow \partial \Omega$ which belongs to the Sobolev space $W^{2,p}\cap \mathcal{C}^{1}(S^{n-1},\RR^{n})$, for $p\geq n-1$. Then
\[h(\Omega)=n-1.\]
\end{thm}

Note that the behavior of entropy for $W^{2,p}$ boundaries for $p<n-1$ is still completely open. From the Sobolev embedding theorem, if $p>n-1$, the assumption on $\mathcal{C}^1$-regularity is superfluous.

\subsection{Outline of the proof and plan of the paper.}

The proofs of theorems \ref{thm:ahlfors} and theorem \ref{thm:sobolev} are completely independent but both have an intense flavor of geometric measure theory.

Sections \ref{sec:prel1} and \ref{sec:proofahlfors} are devoted to Ahlfors regularity and the proof of theorem \ref{thm:ahlfors}, while sections \ref{sec:prel2} and \ref{sec:proofsobolev} deal with Sobolev regularity and the proof of theorem \ref{thm:sobolev}.

The proof of the first main theorem follows two steps. First, we define and compute the volume entropy of some reference Ahlfors $\alpha$-regular convex sets, called Cantor-Lebesgue convex sets and constructed with the familiar Cantor-Lebesgue "the devil's staircase" map. Those convex sets are defined in section \ref{sec:prel1} and the entropy computation is achieved in paragraph \ref{subsec:proofcantorlebesgue}. The second stage of the proof is a comparison argument: taking an arbitrary Ahlfors $\alpha$-regular convex set, we show that its volume entropy is the same as the entropy of the Cantor-Lebesgue convex set of same regularity. Ahlfors $\alpha$-regular maps have a well-defined weak second derivative, which is a measure supported on a Cantor set. In order to compare the volume entropies of the convex domains, we first need to compare the associated Cantor sets. To this end, we use and generalize the main theorem of \cite{saaranen} in paragraph \ref{subsec:comparisoncantor} and we show that the Cantor sets can be related by an order preserving bi-Lipschitz map. The conclusion follows in paragraph \ref{subsec:ahlforsconclusion}.

The second theorem is very close to the main theorem of \cite{bbvc11}, the novelty here is concentrated in lemma \ref{lem:integrable_bound}. Consequently we follow the same outline for the proof. The first step is to show that the (suitably renormalized) Busemann function $\sigma(p)$ converges to the curvature of the boundary as $p$ approaches the boundary. As in \cite{bbvc11}, the main argument lies in the work of D. Alexandroff \cite{Aleksandrov_twice_differentiable_convex}. This convergence is pointwise almost everywhere in \cite{bbvc11} and L$^1$ in this paper (we chose $p\geq n-1$ so that this L$^1$ convergence makes sense). We then apply Lebesgue dominated convergence theorem to deduce that the renormalized volume of balls converge to the so-called centro-projective area (see paragraph \ref{subsec:centroprojectivearea}) which the reader may think of as the total curvature of the boundary. The key lemma \ref{lem:integrable_bound} is specific to the Sobolev regularity context and guarantees the conditions of the Dominated Convergence Theorem.

\subsection{Acknowledgement}

The authors would like to thank Constantin Vernicos, and Marc Troyanov for helpful discussions.

\section{Preliminaries I: Ahlfors regularity and Cantor-Lebesgue domains.}\label{sec:prel1}

\subsection{Ahlfors regularity}\label{subsection:ahlfors}

\begin{defi}
\begin{enumerate}
\item Let $(X,d,\mu)$ be a metric measure space. We say that it is Ahlfors $\alpha$-regular if there exists a constant $C > 0$ such that, for any $x\in X$ and any $r>0$,
\[\frac{1}{C} r^\alpha \leqslant \mu(B(x,r)) \leqslant C r^\alpha. \]
\item Let $\varphi:[0,1]\to \RR$ be a non-decreasing continuous map and let $\mu$ be its derivative in the sense of distributions, seen as a measure. We denote by $X$ the support of this measure. We say that $\varphi$ is Ahlfors $\alpha$-regular if $(X,\abs{\cdot},\mu)$ is Ahlfors $\alpha$-regular.
\end{enumerate}
\end{defi}

\begin{lem}\label{lem:equivalent}
  Let $\varphi:[0,1]\to \RR$ be a continuous function, $\mu$ its derivative and $X$ the support of $\mu$. The function $\varphi$ is Ahlfors $\alpha$-regular if and only if, for any $s,t\in X$, there exists $C>0$ such that,
\[\frac{1}{C}\abs{t-s}^\alpha\leqslant \abs{\varphi(t)-\varphi(s)}\leqslant C\abs{t-s}^\alpha.\]
\end{lem}

This lemma explains why we claimed in the introduction that $\alpha$-Ahlfors regular maps are thought as maps of regularity $\mathcal{C}^{\alpha}$ which are not more regular whenever the second derivative is positive.

\begin{proof}
This follows from the very definition of the derivated measure:
\[\mu(\crochet{s,t})=\varphi(t)-\varphi(s).\]
\end{proof}

Throughout the rest of this section we fix a parameter $p\in \left]2,+\infty\right[$. We now describe the construction and the useful properties of the Cantor-Lebesgue domains $\Omega_p$. In particular we show that the weak second derivative is supported on a $\frac{\log 2}{\log p}$-regular set.

\subsection{Construction of Cantor-Lebesgue domains.}
~

The Cantor-Lebesgue map is constructed inductively.  We start with the identity function $f_{0}:x\mapsto x$, then inductively define
\[f_{n}:[0, 1]\to [0, 1]\]
  \[
  x\mapsto \begin{cases} \frac{ f_{n-1}(p x)}{2} &x < \frac{1}{p}\\
    \frac{1}{2}& \frac{1}{p}\leq x < \frac{p-1}{p}\\
    \frac{f_{n-1}(p x-1)}{2} + \frac{1}{2}&\frac{p-1}{p}\leq x.
  \end{cases}
\]

\begin{figure}[!h]
    \centering
    \begin{subfigure}[b]{0.3\textwidth}
        \includegraphics[width=\textwidth]{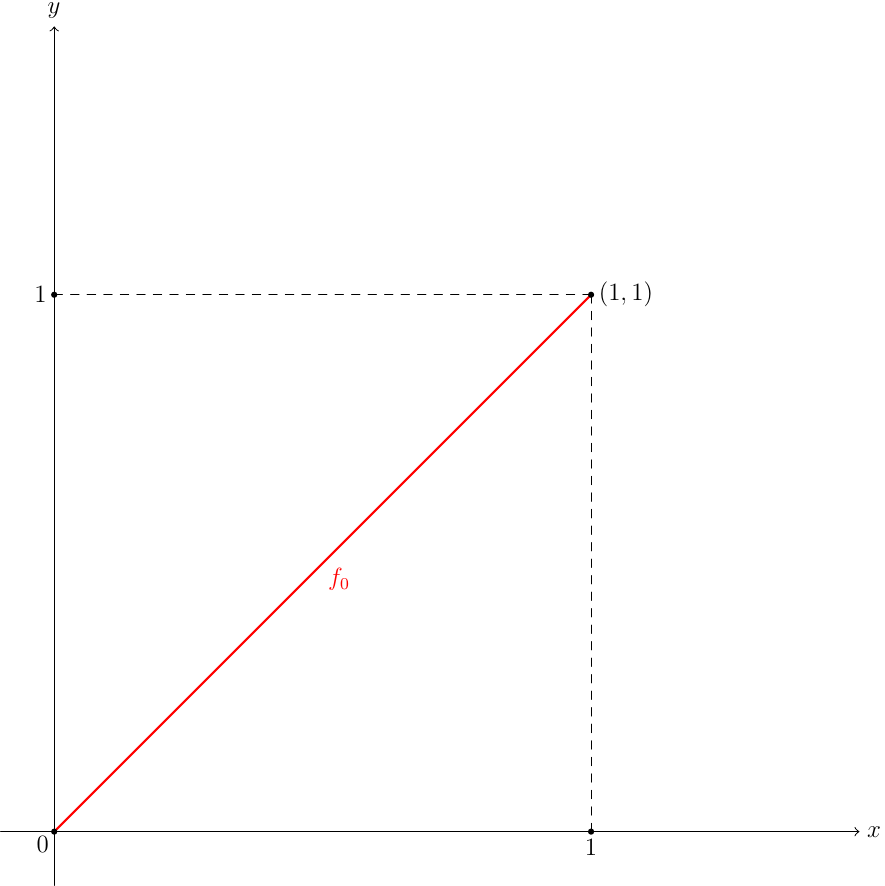}
        \caption{Step 0}
    \end{subfigure}
    ~ 
    \begin{subfigure}[b]{0.3\textwidth}
        \includegraphics[width=\textwidth]{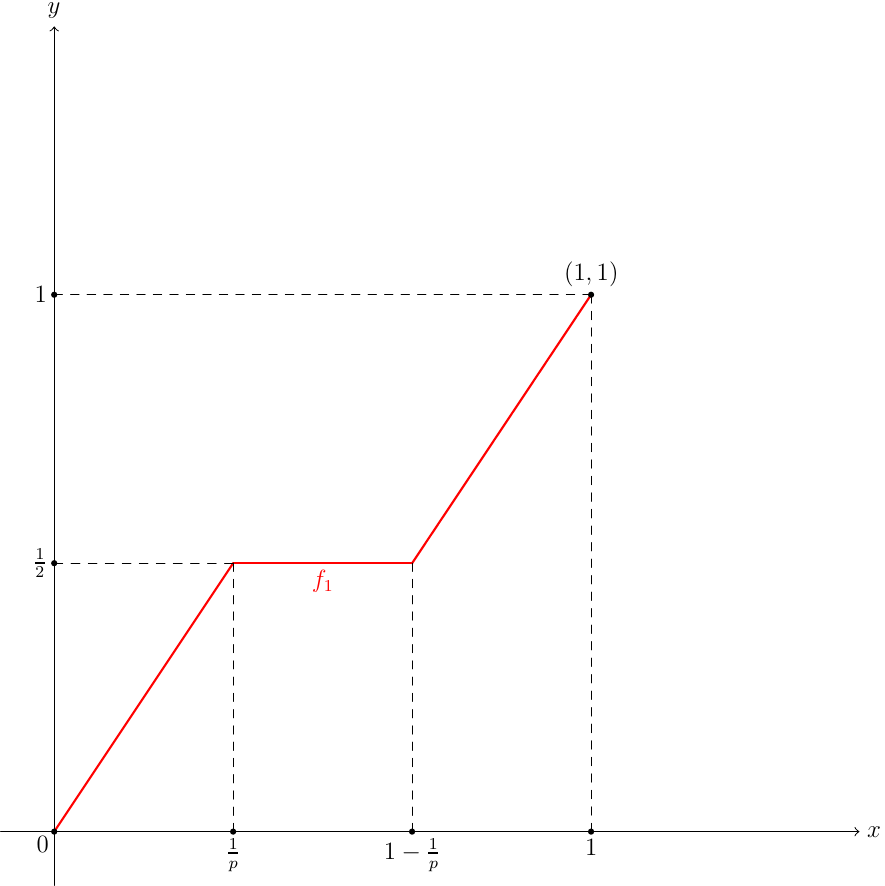}
        \caption{Step 1}
    \end{subfigure}
    ~ 
    \begin{subfigure}[b]{0.3\textwidth}
        \includegraphics[width=\textwidth]{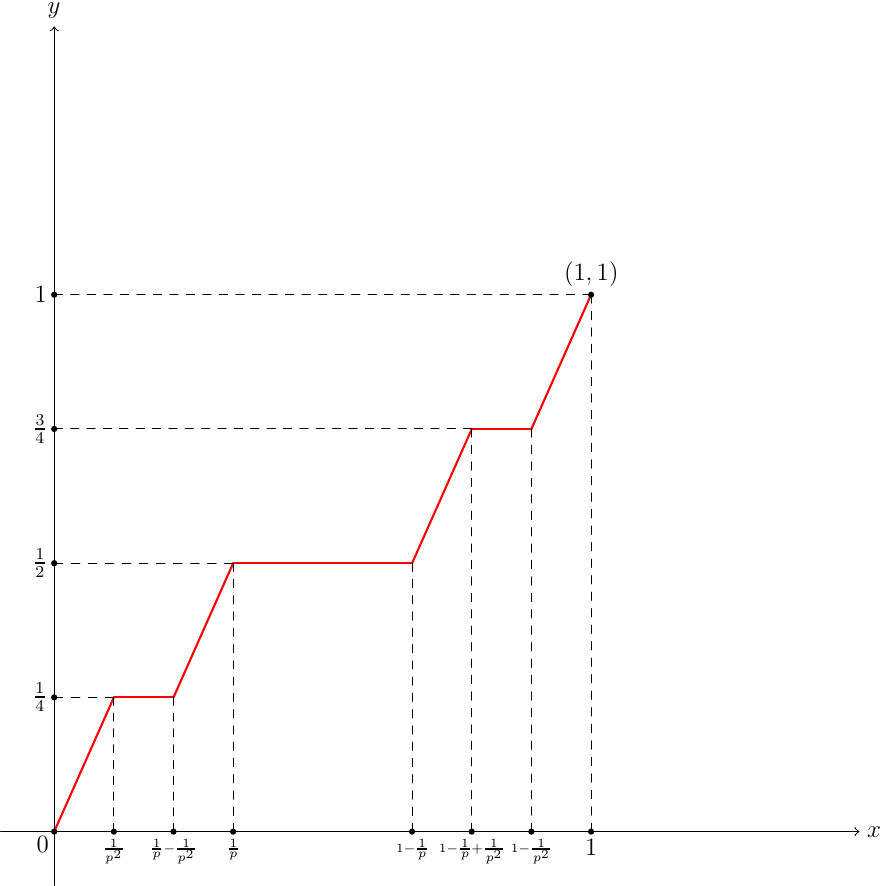}
        \caption{Step 2}
    \end{subfigure}
    \caption{The sequence $(f_n)$}
\end{figure}

The following facts are well known:

\begin{lem}
\begin{enumerate}
\item The sequence $(f_n)$ converges uniformly to a non-decreasing continuous map $f$.
\item The limit map $f$ is $\frac{\log 2}{\log p}$-Hölder continuous.
\end{enumerate}
\end{lem}

\proof  Let $\alpha=\log 2/\log p$.  We suppose inductively that $|f_{n}-f_{n-1}|\leq 2^{-n}\frac{p-2}{p}$.  Then by construction $|f_{n+1}-f_{n}|\leq 2^{-1} |f_{n}-f_{n-1}|\leq 2^{n+1}\frac{p-2}{p}$. Finally the initial step
\[|f_{1}-f_{0}|\leq \frac{1}{2}- \frac{1}{p} = 2^{-1}\frac{p-2}{p}\] is clearly satisfied. Hence  $(f_{n})$ form a uniform Cauchy sequence of continuous maps, converging to some function $f$.   

Now suppose $p^{-k-1}\leq |x-y|\leq p^{-k}$.  Consider
\begin{align*}
  |f(x)-f(y)|&\leq |f_{k}(x)-f_{k}(y)|+2\cdot 2^{1-k}\frac{p-2}{p}.  
\end{align*}
The biggest gap of two points separated by distance $p^{-k}$ for $f_{k}$ is $2^{-k}$. Consequently
\[  |f_{k}(x)-f_{k}(y)|\leq 2^{-k}.\]
And finally
\[|f(x)-f(y)|\leq 5 \cdot 2^{-k}\leq 10 p^{-\alpha(k+1)}\leq 10 |x-y|^{\alpha}.\]
\endproof

\subsection{A probability measure on the Cantor set.}
~

The function $f$ is continous and, as such, it has a weak derivative. It is a Borel measure $\mu$ given by
\[\mu(\crochet{a,b})=f(b)-f(a),\]
Since $(f_n)$ converges uniformly, we can also view $\mu$ as the  limit of the sequence $(\mu_n)$ of derivatives of the $f_n$'s. The supports of the measures $\mu_n$ are decreasing. Then
\[\mathrm{Supp}(\mu)=\bigcap_n\mathrm{Supp}(\mu_n)\]
is a Cantor set. We denote this Cantor set by $C_p$.

\begin{lem}
The measure $\mu$ is $\frac{\log 2}{\log p}$-Ahlfors regular.
\end{lem}
\proof 
By virtue of being H\"older continuous, we have the upper bound for Ahlfors regularity.  To prove the lower bound, let $x\in C_{p}$; it can be written as the countable intersection of closed intervals $\{x\}=\bigcap_{n}I_{n,k_{n}}$ of size $p^{-n}$.  By construction the measure $\mu$ applied to any Cantor interval of size $p^{-n}$ yields $2^{-n}$.  Hence for $2p^{-n}\leq r<2 p^{1-n}$ we have
\[\mu(B(x,r))\geq 2^{-n}\geq  (2p)^{-\alpha}r^{\alpha}.\]
With lemma \ref{lem:equivalent} above, this also shows that $f$ is $\frac{\log 2}{\log p}$-Ahlfors regular.
\endproof
%
%
%
%
%
%
%
%
Now we can explicitly construct the domains. Consider the map $e^{i\pi f}$. Its values give unit vectors in $\RR^{2}$.  If we define $\gamma(t)=\int_{0}^{t}e^{i\pi f}$, this a curve.  If we take a second copy and rotate by $\pi$ around the point $(\gamma(0)+\gamma(1))/2$, we obtain the boundary of a set, which is convex by the monotonicity of $f$. This explicitly defines a map $F:\RR/\mathbb{Z} \cong S^{1}\to \CC\cong \RR^{2}$ by
\[f: t\mapsto \begin{cases}
    \int_{0}^{t}\exp(i\pi f(2 \tau))\:d\tau & 0\leq t < 1/2\\
    \int_{0}^{1/2}\exp(i\pi f(2\tau))\: d\tau - \int_{0}^{t-1/2}exp(i \pi f(2(\tau)\:d\tau & 1/2\leq t < 1
\end{cases}\]
and extending periodically.

This set is the \emph{Cantor-Lebesgue domain}. Throughout the rest of this text, we denote this set by $\Omega_p$.

\begin{defi}\label{def:convexahlfors}
Let $\Omega$ be a convex relatively compact set of $\mathbb{R}^2$. We say that $\Omega$  (or sometimes $\partial\Omega$) is Ahlfors $\alpha$-regular if $\partial\Omega$ could be written locally as the graph of a $\mathcal{C}^1$ function $\varphi :\mathbb{R}\rightarrow\mathbb{R}$ whose derivative is an Ahlfors $\alpha$-regular map.
\end{defi}

Note that, in this case, the curvature measure of $\Omega$ (see the definition at the beggining of paragraph \ref{subsec:centroprojectivearea}) is a Ahlfors $\alpha$-regular measure supported on a Cantor set of dimension $\alpha$. Be aware that the definition of Ahlfors $\alpha$-regularity refers to the weak second derivative and not the first (the notation Ahlfors $1+\alpha$-regular seems too much misleading).

\section{Entropy of Hilbert geometries whose boundaries are Ahlfors regular}\label{sec:proofahlfors}

\subsection{Entropy of the Cantor-Lebesgue domains}\label{subsec:proofcantorlebesgue}

\begin{thm}\label{entropystandard}
For every $p\in ]2,\infty[$ the Cantor Lebesgue domain $\Omega_{p}$ has entropy
\[\ent(\Omega_{p})=\frac{2\alpha}{1+\alpha}\]
where $\alpha=\log 2/\log p$
\end{thm}

\begin{rmk}\label{length}
  By theorem 2.14 in \cite{bbvc11}, as a corollary of the co-area inequalities, the authors show that the entropy is also given by the asymptotic exponential growth rate of volume of \textit{spheres}. Hence, from now on, we focus of computing the length of the circles of radius $R$ and the exponential growth of this length will give the entropy. In fact by Lemma 3 in \cite{entropyhilbertc2} we can consider the length of the scaled boundary $\eta\D \Omega$ ($\eta<1$), which we do.
\end{rmk}

%
\begin{lem}
  Let $\Omega\subset \RR^{2}$ be a convex domain containing the origin.   Then 
  \[\H_{1}(\eta\partial\Omega)=\frac{1}{2} \eta\int_{S^{1}}\frac{1}{|\eta\gamma(\theta)-F^{+}(\theta,\eta)|}+\frac{1}{|\eta\gamma(\theta)-F^{-}(\theta,\eta)|}|d\gamma|\:d\theta.\]
  where $F^{+}(\theta,\eta)$ is the point on $ \D \Omega$ given by the intersection of the tangent line at $\eta\gamma(\theta)$ of $\eta \D \Omega$ and $\D \Omega$, in the positive orientation, and $F^{-}$ is that in the negative orientation. 
\end{lem}
\proof
This is just a direct calculation of the length with the Hilbert norm of the curve
\[\theta\mapsto \eta\gamma(\theta),\]
which is given by
\[\int_{S^{1}}\|\dot{\gamma}\|_{\theta}\:d\theta=\int_{S^{1}}\frac{1}{2}\left( \frac{1}{|\eta\gamma(\theta)-F^{+}(\theta,\eta)|}+\frac{1}{|\eta\gamma(\theta)-F^{-}(\theta,\eta)|} \right)\eta|\dot{\gamma}(\theta)|\:d\theta.\]
\endproof
%
We will assume that our map $\gamma:S^{1}\to \D \Omega$ has unit speed.  Let $\varphi(\theta)$ denote the generating function; so that:
\[\gamma(\theta)=\int_0^\theta e^{i\varphi(s)}ds.\]

\begin{lem}
  Let $\Omega\subset \RR^{2}$ be a convex domain whose boundary generating function $\varphi$ is given by the cumulative distribution function of an Ahlfors $\alpha$-regular measure for some $0<\alpha<1$, and  let $X$ denote the support of this measure, then there is a number c such that  for every $\theta$
\begin{multline*}
  c^{-1}\frac{1+o(1)}{\dist_{+}(\theta,X)+(1-\eta)^{1/(\alpha+1)}}\\\leq |F^+(\theta,\eta)-\eta\gamma(\theta)|^{-1}\\ \leq c \frac{1+o(1)}{\dist_{+}(\theta,X)+(1-\eta)^{1/(\alpha+1)}}
\end{multline*}
as $\eta\to 1$.
\end{lem}
\proof
For each $\theta$ we want to find the point forward for which $\<\eta \gamma(\theta), \nu(\gamma(\theta))\>=\<\gamma(\theta+h), \nu(\gamma(\theta))\>$, i.e.
\begin{align*}
  \<(\eta-1)\gamma(\theta),\nu(\gamma(\theta))\>&=\<[\gamma(\theta+h)-\gamma(\theta)], \nu(\gamma(\theta))\>\\
    &=\<e^{i\varphi(\theta)}\int_{0}^{h}e^{i(\varphi(\theta+t)-\varphi(\theta))}\:dt , \nu(\gamma(\theta)\>\\
  &=\int_{0}^{h}\sin(\varphi(\theta+t)-\varphi(\theta))\:dt.
\end{align*}
Now we use the taylor approximation to $\sin$ to yield
\begin{multline*}
  \int_{0}^{h}[\varphi(\theta+t)-\varphi(\theta)]-\frac{1}{6}[\varphi(\theta+t)-\varphi(\theta)]^{3}\:dt\\\leq\int_{0}^{h}\sin(\varphi(\theta+t)-\varphi(\theta))\:dt\\\leq \int_{0}^{h}\varphi(\theta+t)-\varphi(\theta)\:dt.
\end{multline*}
The generating function $\varphi$ is the cumulative distribution function of an Ahlfors $\alpha$-regular measure on $X$, where $\alpha=\log 2/ \log p$.  Consequently, there is a number $c_{1}\geq 0$ such that $\theta\in [0,2\pi]\setminus X$, and $0\leq t\leq 1/2$ we have
\[c_{1}^{-1}\max\{t-\dist_{+}(\theta,X),0\}^{\alpha}\leq\varphi(\theta+t)-\varphi(\theta)\leq c_{1} \max\{t-\dist_{+}(\theta,X),0\}^{\alpha}\]

Now we apply this to the previous bound to get for $h< \dist_{+}(\theta,X)$ that the integral is zero, and for $h\geq \dist_{+}(\theta,X)$
\begin{align*}
 & \int_{0}^{h}\sin(\varphi(\theta+t)-\varphi(\theta))\:dt\\
  &\geqslant \int_{0}^{h-\dist_{+}(\theta,X)}\left(c_{1}^{-1}s^{\alpha}-s^{3\alpha}/6\right)ds\\
  &\geqslant c^{-1}_{1}([h-\dist_{+}(\theta,X)]^{\alpha+1}/(\alpha+1)-[h-\dist_{+}(\theta,X)]^{3\alpha+1}/(3\alpha+1))\\
  &\geqslant c_{2}([h-\dist_{+}(\theta,X)]^{\alpha+1},
\end{align*}
and
\begin{align*}
  \int_{0}^{h}\sin(\varphi(\theta+t)-\varphi(\theta))\:dt\leq c_{1}[h-\dist_{+}(\theta,C)]^{\alpha+1}/(\alpha+1).
\end{align*}
Now if we can bound $h$ such that 
\[\int_{0}^{h}\sin(\varphi(\theta+t)-\varphi(\theta))\:dt=(1-\eta)\gamma(\theta)\cdot\nu(\theta),\]
but by virtue of convexity the scalar product to the normal to $\D \Omega$ with the position is bounded above and below away from $0$: there is a number $c_{3}$ such that for every $\theta\in [0,2\pi]$
\[c_{3}^{-1}\leq\gamma(\theta)\cdot \nu(\theta)\leq c_{3}. \]
Consequently
\[(1-\eta)c_{3}^{-1}\leq c_{1}[h-\dist_{+}(\theta,X)]^{\alpha+1}/(\alpha+1),\]
and so
\[(1-\eta)^{1/(\alpha+1)}c_{4}+\dist_{+}(\theta,X)\leq h,\]
and 
\[[h-\dist_{+}(\theta,X)]^{\alpha+1}c_{2}\leq c_{3}(1-\eta),\]
so consequently
\[\dist_{+}(\theta,X)+c_{4}(1-\eta)^{1/\alpha+1}\leq h\leq \dist_{+}(\theta,X)+c_{5}(1-\eta)^{1/(\alpha+1)}.\]

Now we must estimate 
\[|\eta\gamma(\theta)-\gamma(\theta+h)|=(1-\eta)\gamma\cdot T(\theta)+\int_{0}^{h}\cos(\varphi(\theta+t)-\varphi(\theta))\:dt\]
We can estimate $1/2 \leq\cos(x)\leq 1$ for $x\in [0,1/2]$ and we bound $\gamma(\theta)\cdot T(\theta)\leq 1$.  Then we arrive at
\begin{multline*}
 \frac{1}{2}\dist_{+}(\theta,X)+c_{4}(1-\eta)^{1/(\alpha+1)}] -(1-\eta) \\\leq |\eta \gamma(\theta)-\gamma(\theta+h)|\\\leq \dist_{+}(\theta,X)+c_{5}(1-\eta)^{1/(\alpha+1})+(1-\eta)
\end{multline*}
If we invert everything we get the result.
\endproof

We note that a similar result holds for $F^{-}(\theta,\eta)$.

\proof[Proof of theorem \ref{entropystandard}]  Now we take advantage of the fact that the standard measure on the Cantor-Lebesgue set is Ahlfors $\alpha$-regular for $\alpha=\log(2)/\log(p)$.  We consider a division of the complement of $X$ into a union of intervals of progressively smaller sets.  At generation $N$ the size of the sets is $(p-2)p^{-N}$, and there are $2^{N}$ of them.
Now we note that for an antipodal set $B(0,R)$ is given by $\tanh(R)\partial\Omega_p$. Using $\tanh R=1-e^{-2R}+O(e^{-4R})$ and setting $\eta=\tanh(R)$ yields
\[\frac{1}{|\tanh(R)\gamma(\theta)-\gamma(\theta+h)|}\sim \frac{1+o(1)}{\dist_{+}(\theta,X)+e^{-2R/(\alpha+1)}}.\] 
Because all terms uniformly constant in $R$ will converge to $1$ as we take the $1/R$ power, it is suffiction to consider
\[\int_{0}^{2\pi}\frac{1}{\dist_{+}(\theta,X)+e^{-2R/(\alpha+1)}}.\]
Now we can break this sum into the generations to yield
\[
  \int_{0}^{2\pi}\frac{1}{\dist_{+}(\theta,X)+e^{-2R/(\alpha+1)}}\:dx =\sum_{N=1}^{\infty} 2^{N}\int_{0}^{(p-2)p^{-N}}\frac{1}{x+e^{-2R/(\alpha+1)}}\:dx\]
\[=\sum_{N=1}^{\infty} 2^{N}\log\left(1 + \frac{(p-2)e^{\frac{2R}{\alpha+1}}}{p^N} \right).\]
In order to proceed we will have to estimate the sum by splitting it in two.  Let
\[\beta=\frac{1}{(\alpha+1)\log p}=\frac{1}{\log 2 +\log p}.\]
We investigate separately the sum over indices $N$ where $N<2R\beta$ and where $N\geq 2R\beta$.
Indeed
\[\sum_{N=1}^{\infty} 2^{N}\log\left(1 + \frac{(p-2)e^{\frac{2R}{\alpha+1}}}{p^N} \right)= \sum_{N=1}^{\infty} 2^{N} \log\left(1+(p-2)p^{2R\beta-N}\right)\]
and the two split sums will have very different asymptotic behaviour.
\begin{align*}
 & \sum_{N=1}^{\infty} 2^{N} \log\left(1+(p-2)p^{2R\beta-N}\right)\\
  &=\sum_{N< 2R\beta}  2^{N} \log\left(1+(p-2)p^{2R\beta-N}\right)+ \sum_{N\geq 2R\beta}  2^{N} \log\left(1+(p-2)p^{2R\beta-N}\right).
\end{align*}
We estimate the first term below coarsely
\begin{align*}
 \sum_{N<2R\beta}2^{N} \log\left(1+(p-2)p^{2R\beta-N}\right) \geq 0
\end{align*}
And from above we estimate coarsely (if $x\geqslant 1$, then $\log(1+x)\leqslant 1+ \log x$):
\begin{align*}
&\sum_{N<2R\beta}2^{N} \log\left(1+(p-2)p^{2R\beta-N}\right)\\
& \leqslant \sum_{N<2R\beta}2^{N} \left(1 + (2R\beta -N)\log p + \log (p-2) \right)\\
& \leqslant \sum_{N<2R\beta}2^{N} \left(1 + 2R\beta\log p + \log (p-2) \right)\\
& \leqslant  2^{\lfloor 2R\beta\rfloor +1}\left(1 + 2R\beta\log p + \log (p-2) \right)
\end{align*}
Now for the second term, we first use the approximation $\log(1+x)\leq x$ for any $x$. We have
\begin{align*}
  \sum_{N\geqslant 2R\beta}2^{N}\log\left(1+(p-2)p^{2R\beta-N}\right) & \leq p^{2R\beta}(p-2)\sum_{N\geqslant 2R\beta}(2/p)^{N}\\
  & \leq  (p-2)(2/p)^{2R\beta}p^{2R\beta}\\
  & =  (p-2) 2^{2R\beta},
\end{align*}
For the lower bound, we want to use
and $x/2\leq \log(1+x)$ for $\abs{x}<1$.  If $2R\beta$ is an integer, the first term of the sum is $\log(1+(p-2))$ and $p-2$ is not necessarily smaller than 1. To avoid this situation we split again the sum into two parts and take off the two first terms.
\begin{eqnarray*}
\sum_{N\geqslant 2R\beta}2^{N}\log\left(1+(p-2)p^{2R\beta-N}\right) & \geqslant & \sum_{N\geqslant 2R\beta +2}2^{N}\log\left(1+(p-2)p^{2R\beta-N}\right)\\
& \geqslant & \frac{(p-2)p^{2R\beta}}{2} \sum_{N\geqslant 2R\beta +2}\left(\frac{2}{p}\right)^{N}\\
& = & \frac{(p-2)p^{2R\beta}}{2}\frac{p}{p-2}\left( \frac{2}{p}\right)^{2R\beta+3}\\
& = & 2^{2R\beta+2}p^{-2}.
\end{eqnarray*}

Consequently we have a two-side bound:
\[\H_{1}(\D B(0, R))\leq (p-2) 2^{2R\beta} + 2^{\lfloor 2R\beta\rfloor +1}\left(1 + 2R\beta\log p + \log (p-2) \right)\]
and
\begin{align*}
  \H_{1}(\D B(0, R))\geq  2^{2R\beta+2}p^{-2}. 
\end{align*}
If we take the $R^{\text{th}}$ root and let $R$ go to infinity we arrive at
\[\lim_{R\to \infty} \left(\H_{1}(\D B(0,R))\right)^{1/R}= 2^{2\beta}=e^{2\log 2/(\log 2+\log p)}=e^{2\alpha/(\alpha+1).}\]
Taking the logarithm yields the result.\endproof

We now come to the general case of a general Ahlfors $\alpha$-regular convex set (recall definition \ref{def:convexahlfors}). The aim of the rest of this section is to prove the first main theorem,

\begin{thm}
Let $\Omega$ be a 2-dimensional Ahlfors $\alpha$-regular convex set. Then its volume entropy satisfies
\[h(\Omega)=\frac{2\alpha}{\alpha+1}.\]
\end{thm}

\subsection{Comparison of Ahlfors regular sets}\label{subsec:comparisoncantor}

In this section, we compare Ahlfors regular Cantor sets. Let us start by fixing notations. We denote by $\kappa$ an Ahlfors $s$-regular measure supported on an Ahlfors $s$-regular set $E$. In the next paragraph we will apply our result below where $\kappa$ is the curvature measure of a $\alpha$-regular Cantor set.

We denote a standard Cantor set (instead of Cantor-Lebesgue Cantor set) by $C_t$. By definition, such a Cantor set is obtained as the support of the derivative of a Cantor-Lebesgue function which is Ahlfors $t$-regular.

\begin{thm}[ordered $s$-regular embeddings]\label{thm:ordered_embedding}
Let $E$ be a totally ordered $s$ regular set.  Let $F$ be a totally ordered $t$ regular set. Assume one of the two regular sets is a standard Cantor set. If $s<t$, there is an order preserving bi-Lipschitz embedding $\varphi:E\to F$ whose constants depend only on the diameters of $E$ and $F$, their Ahlfors constants, and $t$ and $s$.
\end{thm}

\begin{rmk}
This statement is an extension of the theorem 3.3 in \cite{saaranen}. Here we add the fact that we can choose the bi-Lipschitz map $\varphi$ to be order-preserving. Note that our proof uses the result of \cite{saaranen}.
\end{rmk}

Starting with the map $f:E\to F$ given by theorem 3.3 of \cite{saaranen}, we proceed as follows
\begin{enumerate}
\item We first construct binary trees for which $E$ and $F$ are topological boundaries of the trees. The construction is inspired by a paper of F. Choucroun \cite{choucroun} but is not quite the same.
\item We choose an (incomplete) metric on the tree so that the associated metric (completion) on the boundary is bi-Lipschitz equivalent to the original metric on the Cantor set, given as a subset of $\mathbb{R}$. The fact that every Cantor set is bi-Lipschitz equivalent to the metric boundary of some \textit{binary} tree seems new.
\item We then extend the given map $f$ by a bi-Lipschitz map $\tilde{f}$ between the trees.
\item We finally compose $f$ with an automorphism of the standard tree to reorder the boundary.
\end{enumerate}

\textbf{Step one : construction of the binary tree}
~

Let $E$ be an ordered Ahlfors $\alpha$-regular set. As being closed, we can write the complement of $E$ as a disjoint union of open intervals:
\[E^c=\bigcup_{i=1}^\infty I_i.\]
On the $I_i=\crochet{a_i,b_i}$, we give the order such that $b_i-a_i\geqslant b_{i+1}-a_{i+1}$ and, if $b_i-a_i = b_{i+1}-a_{i+1}$, then $b_i<a_{i+1}$.


Set $A_0=\crochet{0,1}$ and, for $k>0$,
\[A_k=\crochet{0,1}\setminus\bigcup_{j=1}^k I_j.\]
We note that $A_k$ is a disjoint untion of $k+1$ closed intervals, denoted $v_0^k,\cdots,v_k^k$.
Let us consider
\[V_n=\bigcup_{k=0}^{n}\set{v_0^k,\cdots,v_k^k}.\]
We build a finite tree for which the vertices are the elements of $V_n$. We now describe the set of edges. We first remark that the union defining $V_n$ is not disjoint: passing from $A_k$ to $A_{k+1}$ creates a hole in one of the $v_i^k$. So precisely there are in $\set{v_0^{k+1},\cdots,v_{k+1}^{k+1}}$ exactly 2 new elements and there is precisely one element in $\set{v_0^k,\cdots,v_k^k}$ not in $\{v_{0}^{k+1},\cdots, v_{k+1}^{k+1}\}$. Hence there exist $\alpha_k\in \set{v_0^k,\cdots,v_k^k}$ and $\beta_k^1,\beta_k^2 \in \set{v_0^{k+1},\cdots,v_{k+1}^{k+1}}$ such that
\[\set{v_0^{k+1},\cdots,v_{k+1}^{k+1}} = \set{v_0^k,\cdots,v_k^k} \setminus \set{\alpha_k}\cup \set{\beta_k^1,\beta_k^2}\] and the union is disjoint. For each $k$, we place two edges: $\crochet{\alpha_k,\beta_k^1}$ and $\crochet{\alpha_k,\beta_k^2}$. We get a finite binary tree $T_n$.
Finally we consider the $\mathbb{R}$-tree $T=\bigcup_{n=0}^\infty T_n$. If needed we put the Cantor set in the notation and denote the tree by $T_E$.\\

\textbf{Claim} Every vertex has exactly two descendants.

\begin{proof}
From the construction, it is obvious that each vertex has either 0 or 2 descendants. Suppose there exists $v_i^k$ with no descendant. Then $v_i^k\subset E$ and $E$ is not totally discontinuous; a contradiction.
\end{proof}
Let us now endow the binary tree $T$ with a metric.  Given a vertices $v_{i}^{k}$ and $v^{l}_{m}$, we define the distance between them to be the size $|v|$ of the smallest closed interval which contains both of them. This metric is incomplete, the vertices accumulate on the boundary which is at finite distance (one) to the basepoint $v_0^0$. We denote by $d_T$ the metric on $\partial T$ induced by the completion of the metric described above.\\

\textbf{Step two : comparison of metrics on the regular set}
~

\begin{thm}
  The tree metric $d_{T}$ and Euclidean metric $d$ are bi-Lipschitz equivalent on $E$.
\end{thm}
\proof The proof will proceed in several parts.  First we will show that Ahflors regularity implies that for a given vertex $v$ in our tree with descendents $u$ and $w$, it follows that there is a number $K>0$ which depends only on $s$ and the proportionality constants of $s$-regularity, such that
\[\max\{|u|,|w|,|v|-|u|-|w|\}\leq K\min\{|u|,|w|,|v|-|u|-|w|\}.\]
This requires judicious application of the Ahlfors regularity, and a Vitali covering argument. 

First we will prove that every closed interval $u$ contains a ball of radius $|u|/4$.  First suppose the open interval $I$ contained in $u$ is bigger than $|u|/3$.  It follows that $I$ is smaller or equal in size to any other open interval, so the open intervals on either side of $u$ are bigger in length than $|u|/3$.  Now $u=v\cup I\cup w$.  Without loss of generality assume that $|v|\geq |w|$.  Then let $v=[a,b]$  so that $B(b,|u|/3)\cap E\subset u\cap E$.  Now suppose $|I|<|u|/3$. Without loss of generality we may assume $|v|\geq |u|/3$, and so if $v=[a,b]$ then $B(b,|u|/3)\cap E\subset u\cap E$.

As a result we can deduce that
\[c|u|^{s}/3^{s}\leq\kappa(u)\leq C |u|^{s}\]
(recall that $\kappa$ denotes the curvature measure on the Cantor set).
Now consider some closed interval $u=[a,b]$ with $I_{k}$ on its left and $I_{l}$ on its right. Without loss of generality assume that $|I_{l}|\leq |I_{k}|$.  We have two cases, either $|I_{l}|\leq |u|$, or $|I_{l}|>|u|$.  If we have the latter then  $B(b,|I_{k}|)\cap E=u\cap E$, and
\begin{align*}
0&= \kappa(B(b,|I_{k}))-\kappa(u)\\
&\geq c |I_{k}|^{s}-C |u|^{s},
\end{align*}
which implies 
\[\frac{|I_{k}|}{|u|}\leq \left(\frac{C}{c}\right)^{1/s}.\]

Next consider some interval $u=v\cap I\cap w$.  Then $\kappa(u)=\kappa(v)+\kappa(w)$ and
\[\kappa(u)\geq c(|v|+|I|+|w|)^{s}/3^{s},\]
while 
\[\kappa(u)=\kappa(v)+\kappa(w)\leq C(|v|^{s}+|w|^{s})\leq C'(|v|+|w|)^{s}.\]
From this it follows that
\[|I|\leq K (|v|+|w|),\]
where $K$ dependts only on $s$, $c$ and $C$.

Finally cover $u$ in balls of radius $|I|$ and centered in $E\cap u$.  This is possible, because $|I|$ is the largest gap in $u\cap E$.  First we have a finite sub-cover.  We can then take a Vitali cover to yield balls $B(x_{i},|I|)$, $i=1,\ldots n$, of disjoint balls, such that $\{B(x_{i},3|I|):\:i=1,\ldots n\}$ covers $u$.  Consequently $n\geq |u|/6|I|.$  But then
\begin{align*}
  \kappa(u)&\geq \sum_{i}\kappa(B(x_{i},|I|))\\
  &\geq n c |I|^{s}/3^{s}\\
  &\geq  3^{-(1+s)} c |u||I|^{s-1}.
\end{align*}

Similarly
\[\kappa(u)\leq C|u|^{s}.\]
From which it follows that
\[\frac{|I|}{|u|}\geq K,\]
where $K$ and $C'$ depend only on $s$,$c$ and $C$.

It is now clear how we can show that $d$ and $d_T$ are bi-Lipschitz equivalent. Indeed, let $x,y\in E$ and let $v$ be the smallest closed interval containing both ot $x$ and $y$ and $w$ be the biggest open interval contained in $\crochet{x,y}$. We have
\[\frac{1}{K} d_T(x,y)\leqslant \abs{w}\leqslant \abs{x-y}\leqslant \abs{v} = d_T(x,y).\]

\endproof
\textbf{Step three : extension of $f$}
~

We are guaranteed a bi-Lipschitz map from $E$ to $C_{t}$ for any standard Cantor set $t>s$ by Theorem 3.3  in \cite{saaranen}.  Let $f$ denote this function.  Let $\tilde{f}:T_{E} \to T_{C_t}$ be given by 
\[v\mapsto \cap \{u\in T_{C_t}:u\supset f(v)\}.\]
Let $I\subset \tilde{f}(v)$ be the largest open interval (with possible semantic order) contained in $\tilde{f}(v)$.  By construction there are an $x,y \in v\cap E$ such that $f(x)$ and $ f(y)$ are separated by $I$, so
\[|I|\leq |f(x)-f(y)|\leq C |x-y|,\]
where $C$ is the bi-Lipschitz constant of the map.  But there is a number $K$ which depends only on $t$ such that $|\tilde{f}(v)|\leq K |I|$.  
Now given two vertices $u'$ and $u''$ in $T_{E}$.  Let $J$ denote the maximal open set separating them. Then $|J|\leq d_{T}(u',u'') \leq |v|\leq K |J|$.  Now we know that there is an $x\in u'$ and $ y\in u''$ such that $d(x,y)\geq |J|$, and hence $C|f(x)-f(y)|\geq|x-y|\geq |J|$.  Consequently 
\[d_{T}(\tilde{f}(u'),\tilde{f}(u''))\geq |f(x)-f(y)|\geq C^{-1}|J|\geq C^{-1}K^{-1}d_{T}(u',u'').\]  For the reverse inequality, note that we have a $K'$ depending only on $t$, such that $\tilde{f}(u')$ and $\tilde{f}(u'')$ are separated by a set of size at least $(K{'})^{-1} d_{T}(\tilde{f}(u'),\tilde{f}(u''))$.  Consequently there is an $x\in u'$ and a $y$ in $u''$ such that
\[\frac{1}{K'}d_{T}(\tilde{f}(u'),\tilde{f}(u''))\leq|f(x)-f(y)|\leq C|x-y|\leq C d_{T}(u',u'').\]

\textbf{Step four : reordering}
~

Subsequently, the most useful property of the standard Cantor set is its self similarity.  In particular for the standard Cantor tree $T_{C_t}$ with the ultrametric $d_{T}$, if we flip any branch, it is an isometry of the tree $T_{C_t}$, as is any uniform limit of isometries.  Because our tree is binary we can identify each vertex with a word composed of the letters $l$ and $r$, for left and right.  Now  given $u$ and $v$ which are not descendents we say that $u$ is to the left of $v$ if their minimal closed cover is $w$ is such that $u$ is in the left branch from $w$ and $v$ is in the right branch. We proceed inductively on word length.  Suppose $\phi_{n}\circ\tilde{f}$ preserves the order of the first vertices with word length less than or equal to $n$. For each vertex $w$ of word length $n$ consider its two descendents $u$ and $v$ such that $u$ is to the left of $v$.  If $\phi_{n}\circ\tilde{f}(u)$ is to the left of  $\phi_{n}\circ\tilde{f}(v)$ proceed if not then postcompose $\phi_{n}$ with the branch flip, which flips the tree at the common vertex of $\phi_{n}\circ\tilde{f}(u)$ and $\phi_{n}\circ\tilde{f}(v)$.  This flip leaves every other descendent pair unchanged because this common vertex is a descendent of $\phi_{n}\circ\tilde{f}(w)$.  Let $\phi_{n+1}$ be the end result.

In this way we construct a sequence of isometries which converges uniformly to a limit isometry $\phi$.  Then $\hat{f}=\phi\circ \tilde{f}$ is a bi-Lipschitz map which preserves the semantic order. In fact we get the following
\begin{cor}
  \label{cor:Lipschitz_extension}
  Let $E$ and $E'$ be $s$ and $t$ regular subsets of $[0,1]$.  Suppose $\{0,1\}\subset E$.  There is an order-preserving map $f:[0,1]\to [0,1]$ which is bi-Lipschitz on its image, and such that $f(E)\subset E'$.
\end{cor}
\proof We take a standard Cantor set $C$ with $s<\tau <t$.  Then there is an order preserving bi-Lipschitz map $f:E\to C$ and $g:C\to E'$. First we extend $f$ to $[0,1]$.  We do this by mapping each open  interval $I$ on the complement to the image of it's endpoints and scaling linearly.  Because the map $f$ is bi-Lipschitz and order preserving, so is the extended map.  By composing the extensions we get the desired map.
\endproof

\subsection{Entropy of a domain with curvature an Ahlfors regular measure.}\label{subsec:ahlforsconclusion}

Recall that we dispose of a 2-dimensional convex set $\Omega$, which is Ahlfors $\alpha$-regular: its curvature measure is supported on $E$, an $\alpha$-regular Cantor set. We choose $t>\alpha$ and we let $C_{t}$ denote the standard Cantor set.

We can express the complement of the support of our measure as a union of open intervals 
\[E^{c}=\bigcup_{j}I_{j},\]
where $I_{j}=]a_{j},b_{j}[$, with $a_{j},b_{j}\in E$.  We know from the previous paragraph that, associated to $C_t$, we have a bi-Lipschitz order preserving map $\varphi : E \rightarrow C_t$.

As in the case of standard Cantor set, the lengths of the Finslerian circles of radius $R$ in $\Omega$ are given by a sum of integrals of the form 
\begin{eqnarray*}
\H_{1}(\D B(0, R)) &\sim & \sum_j \int_{a_{j}}^{b_{j}}\frac{1}{(x-b_{j})+e^{-2R/(\alpha+1)}}\:dx\\
& = & \sum_j \int_{0}^{[b_{j}-a_{j}]}\frac{1}{x+e^{-2R/(\alpha+1)}}\:dx.\\
& = & \sum_j \log (( b_{j}-a_{j})e^{2R/(\alpha+1)} +1).
\end{eqnarray*}
But $|b_{j}-a_{j}|\leq C|\varphi(b_{j})-\varphi(a_{j})|$. Consequently we can estimate 
\begin{eqnarray*}
\log(e^{2R/(\alpha+1)}|b_{j}-a_{j}|+1) & \leq & \log(e^{2R/(\alpha+1)}C|\varphi(b_{j})-\varphi(a_{j})|+1)\\
& \leq & \sum_{J\subset ]\varphi(a_{j}),\varphi(b_{j})[}\log(C|J|e^{2R/(\alpha+1)}+1)
\end{eqnarray*}
where $I$ are the open intervals whose disjoint union makes the complement of the $t$ dimensional cantor set. This last inequality follows from the elementary identity
\[\log(1+\lambda +\mu)\leqslant \log(1+\lambda) + \log(1+\mu)\]
whenever $\lambda$ and $\mu$ are positive.
Because our map is order preserving, we can just estimate from above (there will be no open intervals appearing more than once in the image)
\begin{eqnarray*}
\sum_{J\subset E^{c}}\log(1+e^{2R/(\alpha+1)}|J|) & \leq & \sum_{J\subset C_{t}^{c}}\log(1+C e^{2R/(\alpha+1)}|J|)\\
 & = & \sum_{n=1}^\infty \sum_{\abs{J}=\frac{p-2}{p^n}}\log(1+C(p-2) e^{2R/(\alpha+1)}p^{-n})\\
 & = & \sum_{n=1}^\infty 2^n \log(1+C(p-2) e^{2R/(\alpha+1)}p^{-n})
\end{eqnarray*}
Here $p$ is the real number such that $\frac{\log 2}{\log p}=t=\dim C_t$.
Then we proceed as before. Setting $\beta=\frac{1}{(\alpha +1)\log p}$, we break the sum into two parts where either $n<2R\beta$ or $n\geqslant 2R\beta$ and use the techniques for bounding above:
\begin{align*}
& \sum_{J\subset C_{t}^{c}}\log(1+C e^{2R/(\alpha+1)}|J|)\\
& \leqslant \sum_{n<2R\beta}2^n\left(1+\log \left(C(p-2)p^{2R\beta -n} \right)\right) + \sum_{n\geqslant 2R\beta}2^nC(p-2)p^{2R\beta - n}
\end{align*}
using for the first term the identity $\log (1+x)\leqslant 1 +\log x$ and for the second $\log(1+x)\leqslant x$.
The first term ($n<2R\beta$) is again
\[\sum_{n<2R\beta}2^n +2^n \log (C(p-2)) + 2^n(2R\beta-n)\log p.\]
which we can coarsely bound by
\[2^{\lfloor 2R\beta\rfloor +1}2R\beta\log p.\]
For the second term ($n\geqslant 2R\beta)$, we have
\[\sum_{n\leqslant 2R\beta}\left(\frac{2}{p}\right)^nC(p-2)p^{2R\beta}\leqslant \left(\frac{2}{p}\right)^{\lfloor 2R\beta\rfloor +1} C(p-2)p^{2R\beta}\leqslant 2^{\lfloor 2R\beta\rfloor +1}C(p-2)\]
Consolidating terms we get that the integral is bounded by 
\[2^{\lfloor 2R\beta\rfloor +1}2R\beta\log p + 2^{\lfloor 2R\beta\rfloor +1}C(p-2)\]
where $t=\log2/\log p$ and $K$ is independent of $R$. Taking the power of $1/R$, the logarithm and letting $R$ go to infinity yields the entropy bound
\[\mathrm{Ent}(E)\leq 2t/(\alpha+1).\]
But $t$ is arbitrarily close to $\alpha$.
We obtain the lower bound in the very same way: we use this time an order preserving map $\psi : C_s \rightarrow E$ and produce a lower bound with the same techniques as for the standard Cantor set.

\section{Preliminaries II: Sobolev regularity and Busemann functions}\label{sec:prel2}

\subsection{What is a Sobolev-regular convex set ?}\label{prel2:parametrization}

\begin{con}
Throughout this text, for a convex set $\Omega$, and a point $x\in\partial \Omega$, $\nu(x)$ denotes the unit inner normal of $\Omega$ and based at $x$.
\end{con}

We clarify the definition of $W^{2,p}$ convex domains $\Omega$ (via charts) and show that the definition is equivalent to requiring that the boundary can be parametrized by a map \[\varphi_{\Omega}:S^{n-1}\to \D \Omega,\]
of class $W^{2,p}\cap C^{1}(S^{n-1},\RR^{n})$, which we define to be the space of continuously differentiable functions $f:S^{n-1}\to \D\Omega\subset \RR^{n}$ whose restriction to coordinate charts, $\psi:U\subset S^{n-1}\to \RR^{n-1}$, denoted by $\varphi_\Omega\circ\psi^{-1}$, are in 
\[W^{2,p}(\psi(U)\subset \RR^{n-1},\RR^{n}).\]

We say a domain $\Omega$ has regularity $W^{2,p}\cap C^{1}$ for $p\geq (n-1)$ if for every $x\in \D \Omega$ there is an $n-1$ dimensional supporting plane $P$ through $x$, and open subsets $U'$ of $P$ and $U$  of $\D \Omega$ containing $x$ such that $U$ is given by the graph of a function $f\in W^{2,p}\cap C^{1}(U')$ i.e. for every $x$ there is a $U'\subset T_{x}\D\Omega$ containing $0$, and a function $h:U'\to \RR^{+}\cup\{0\}$ such that the set $ \{z-h(z)\nu(x)+x|\; z\in U'\} $ is a subset of $\D\Omega$ and contains an open neighborhood of $x$ in $\D\Omega$.

Note that for $p=n-1$, $W^{2,p}(S^{n-1})$ does not embed in $C^{1}$ which is why we study parametrizations in $W^{2,p}\cap C^{1}$. If $p>n-1$, $W^{2,p}(S^{n-1})$ does embed in $C^{1}$.

We define the \emph{barycenter} of a domain $\Omega$ as the center of mass of the boundary 
\[\bc(\Omega)=\int_{\D \Omega}x\:d\H_{n-1}(x)/\H_{n-1}(\D \Omega).\]
Projecting the unit sphere around the barycenter onto $\D \Omega$ gives a map.
\[\varphi_{\Omega}:S^{n-1}\to \D \Omega\]
  which takes $\theta\in S^{n-1}$ to the point on $\D \Omega$ in the direction $\theta$ from $\bc(\Omega)$
If $\Omega$ is a bounded convex domain then $\bc(\Omega)$ exists as does $\varphi_{\Omega}$.

\begin{lem}
  \label{lem:sobolev_reg}
  Let $p\geq n-1$.  The bounded convex domain $\Omega$ is $C^{1}\cap W^{2,p}$ regular if and only if the map 
  \[\varphi_{\Omega}:S^{n-1}\to \D \Omega,\]
  is in the space $W^{2,p}\cap C^{1}(S^{n-1},\RR^{n})$.
\end{lem}
\proof
The map $\varphi$ can be identified with the map $\rho:S^{n-1}\to \RR^{+}$ for which
\[\varphi(\theta)=\rho(\theta) \theta+\bc(\Omega).\]
Then $\varphi\in W^{2,p}(S^{n-1},\RR^{n})$ is equivalent to $\rho\in W^{2,p}(S^{n-1})$.

Let $h:U\subset T_{x}\D \Omega\to \RR$ be such that $\{x+z-\nu(x)h(z)\mid z\in U\}\subset \D \Omega$.  Let $x_{0}$ denote the barycenter of $\Omega$.  Then by assumption $h$ is $W^{2,p}\cap C^{1}(U)$.  Consider the map 

\[\Phi:U\subset T_x\D\Omega\to S^{n-1}\subset \RR^{n}\qquad\Phi:z\mapsto \frac{z+x-h(z)\nu(x)-x_{0}}{|z+x-h(z)\nu(x)-x_{0}|},\]
which takes $z$ to $\D \Omega$ along $\nu(x)$, and then radially to the unit sphere centered at $x_{0}$.  This map is $C^{1}$, as it is a quotient of $C^{1}$ functions.  It  is also $W^{2,p}$ because it is a product of $C^{1}\cap W^{2,p}$ functions and $C^{1}\cap W^{2,p}$ is an algebra \cite{Royden_Algebra_Space}.  Let $\zeta=z+x-h(z)\nu(x)$.  Consider the derivative
\begin{align*}
  D\Phi_{z}&=\frac{Id-\nu(x)dh(z)}{|z+x-h(z)\nu(x)-x_{0}|}\\
  &\quad-\frac{(z+x-h(z)\nu(x)-x_{0})}{|z+x-h(z)\nu(x)-x_{0}|}\frac{\<z+x-h(z)\nu(x)-x_{0},Id-\nu(x)dh(z)\>}{|z+x-h(z)\nu(x)-x_{0}|}.\\
\end{align*}
This is the composition of this inverse of projection from $T_{\zeta}\D \Omega\to T_{x}\D \Omega$ and the projection $T_{\zeta}\D \Omega$ to $T_{(\zeta-x_{0})/|\zeta-x_{0}|} S^{n-1}$ but because the line from $x_{0}$ to $ \zeta$  is transverse to $\D \Omega$ at $\zeta$ this projection is invertible.  Consequently by the inverse function theorem $\Phi^{-1}$ is locally $C^{1}$.  Finally let $F=\eta\circ \Phi$ where $\eta:\widehat{U}\subset S^{n-1}\to\RR^{n-1}$ is a smooth chart.  Then we have deduced that $F^{-1}:U\subset \RR^{n-1}\to T_{x}\D \Omega$ is locally $C^{1}$.  Now consider for a $C^{1}$ diffeomorphism between open subsets of $\RR^{n-1}$:
\[D(F^{-1})(x)=(DF)^{-1}(F^{-1}(x)),\]
where $(DF)^{-1}$ is the matrix inverse of $DF$.  This is a composition of continuous functions, so is continuous.  Now taking the derivative with respect to a vector $X$ in $\RR^{n-1}$
\[D_{X}D(F^{-1})(x)=DF^{-1}_{F^{-1}(x)}\cdot(D_{DF^{-1}_{x}X}\cdot DF_{F^{-1}(x)})DF^{-1}_{F^{-1}(x)}\]
where we have used the fact
\[D(A^{-1}(x))=A^{-1}(x)\cdot DA(x)\cdot A^{-1}(x)\]
for any invertible matrix valued function $A:U\subset \RR^{n-1}\to GL(d)$.
The change of variables is at least a $C^{1}$ diffeomorphism and the matrices $DF$ and $DF^{-1}$ are continuous (and hence bounded). This imply that $D_{X}D(F^{-1})(x)$  is bounded uniformly almost everywhere by $\|DF\|_{\infty}\|DF^{-1}\|_{\infty}|D^{2}F|$.  This is valid for any $C^{2}$ approximation of $F$, and so passing to the limit is valid for $F$ \cite[Chap. 5]{PDEs}.   And so if $|D^{2}F|\in L^{p}(T_{x}\D \Omega)$ then $D^{2}(F^{-1})$ is locally in $L^{p}(\RR^{n-1})$.  Finally the map $\varphi$ is the composition of $\Theta:z\mapsto z+x-h(z)\nu(x)$ and $\Phi^{-1}$ which are both $C^{1}\cap W^{2,p}$ i.e. 
\[\varphi=\Theta\circ (\Phi^{-1})\]
on the set $\Phi(U)\subset S^{n-1}$.  Furthermore $\Phi$ is a $C^{1}\cap W^{2,p}$ diffeomorphism.   Hence the composition is $C^{1}\cap W^{2,p}$.

For the reverse implication, let $x$ be a point in $\D\Omega$.  Then consider the map $\Psi=\Pi_{T_{x}\D\Omega}\circ\varphi$ where $\Pi_{T_{x}\D\Omega}$ with orthogonal projection onto the tangent plane at $x$ of $\D\Omega$.  By a similar chain of reasoning $\Psi $ is in $W^{2,p}\cap C^{1}(S^{n-1},T_{x}\D\Omega)$ and is locally invertible around $\varphi^{-1}(x)$. By a similar argument to the forward implication, it too is in $W^{2,p}\cap C^{1}$.  Then $h(z)$ is given by $\<(\varphi\circ \Psi^{-1})(z),\nu(x)\>-\<x,\nu(x)\>$ which is a composition of $W^{2,p}\cap C^{1}$ functions which is $W^{2,p}\cap C^{1}$, \cite{NLglobalanalysis}.
\endproof

\subsection{The curvature measure and the centro-projective area for the Sobolev regularity}\label{subsec:centroprojectivearea}
~

Let $\Omega$ be a convex compact set.
Convex domains naturally have some regularity. They are automatically Lipschitz regular, and twice differentiable almost everywhere but in addition we can define a set valued Gauss map for the boundary of any domain.  For every $x\in \partial\Omega$ and $\nu\in\mathbb{S}^{n-1}$, we consider the halfspace
\[H(x,\nu)= \{y:\<x,\nu\>\leq \<\nu,y\>\}\]
and then we set
\[G(x)=\left\lbrace \nu \in \mathbb{S}^{n-1} \;\mbox{ such that }\; H(x,\nu)\cap \mathrm{int}(\Omega) = \emptyset \right\rbrace.\]
This allows us to define a Gauss curvature measure on $\partial \Omega$ by 
\[\kappa(A)=\mu(\bigcup_{x\in A} G(x)),\]
where $\mu$ is the usual measure on the sphere.
For a $C^{2}$ domain the measure corresponds with the change of variables formula 
\[\kappa(A)=\int_{A}\det DG(x)\:dx.\]
Already the $C^{2}$ assumption implies that the entropy of a domain is $n-1$.  In \cite{bbvc11} this supposition was weakened to $C^{1,1}$.  If one weakens the asumption to $W^{2,p}$ for $p\geq n-1$, we still have a well defined Gauss curvature. Indeed, the determinant of the Gauss map $\det DG$ belongs to L$^1$ because $p\geq n-1>\frac{n-1}{2}$.

Assume now that the origin $o$ of $\mathbb{R}^n$ belongs to $\Omega$ (this is no restriction) and let $a :\partial \Omega \rightarrow \mathbb{R}$ be the positive function such $-a(p)p\in\partial\Omega $ (see the introduction of \cite{bbvc11}). For instance if $o$ is a center of symmetry of $\Omega$, then $a$ is just the constant 1. The letter $a$ stands for antipodal function.
We can now recall the definition of the centro-projective area in the following way.

\begin{defi}[centro-projective area]
Let $\Omega$ be a convex set such that $\partial\Omega$ is parametrized by a map in $W^{2,p}\cap C^{1}(\mathbb{S}^{n-1})$. The centro-projective area is defined as
\[\mathcal{A}(\Omega)=\int_{\partial \Omega}\frac{2a\sqrt{\det DG(x)}}{\left((1+a)\<\nu(x),x\>\right)^{\frac{n-1}{2}}}d\mathcal{L}(x)\]
\end{defi}

\begin{lem}
  Let $\Omega$ be a convex domain of class $C^{1}\cap W^{2,p}$ for $p\geq n-1$.  Then the centro-projective area is nonzero
  \[\A(\Omega)>0.\]
\end{lem}
\proof
The goal is to show that the Gauss map $G$ is of class $W^{1,n-1}\cap C^{0}$. In this case by the change of variables formula 
\cite[Theorem 6.3.2]{GFT} yields
\[\int_{S^{n-1}}\det DG(\theta)\:d\theta=|S^{n-1}|.\]
Thus $\det DG$ is nonzero on a set of positive measure.  Furthermore, by convexity and Alexandrov's theorem on twice differentiability $\det DG$ is non-negative almost everywhere.  Consequently so is $\sqrt{\det DG}$, and hence $\A(\Omega)>0$.

To see that $G$ is in $C^{0}\cap W^{1,n-1}$ let $\psi:\psi^{-1}(U)\subset S^{n-1}\to U\subset \RR^{n-1}$ be a coordinate chart.  Let  $\Phi:U\to \D \Omega$ be the map $x\mapsto \varphi(\psi(x))$. By Lemma \ref{lem:sobolev_reg} it  is of class $C^{1}\cap W^{2,p}$.  Consequently for an orthonormal frame of vector fields $X_{1},\ldots, X_{n}$  the function 
\[D_x\Phi(X_{1})\wedge\cdots \wedge D_{x}\Phi(X_{n-1}),\]
is of class $C^{0}\cap W^{1,n-1}$ as is
\[x\mapsto\star(D_{x}\Phi(X_{1})\wedge\cdots \wedge D_{x}\Phi(X_{n-1}))/|D_{x}\Phi(X_{1})\wedge\cdots\wedge D_{x}\Phi(X_{n-1})|,\]
where $\star$ is the hodge star, because $C^{0}\cap W^{1,n-1}$ is an algebra \cite[Theorem 2.1]{Royden_Algebra_Space}.  But this is the Gauss map.
\endproof

\subsection{Some general facts about Busemann functions}\label{section:busemann}

Let us start by an elementary fact.

\begin{lem}\label{lem:starintegration}
Let $S$ be a bounded star-shaped domain of $\mathbb{R}^n$ with $0$ as a star and with a $\mathcal{C}^1$ boundary. Hence there is a parametrization of $S$ given by
\[\fonction{\Phi}{B(0,1)}{S}{x}{x\varphi(\frac{x}{\abs{x}})}\]
where $\varphi$ is a real valued function on $\mathbb{S}^{n-1}$.
Then
\[\Vol (S)= \int_{\mathbb{S}^{n-1}}\varphi(x)^n dx\] 
\end{lem}

This will be very useful for computing the Euclidean area of a Finslerian ball (hence evaluationg the Busemann function).

For the rest of this paragraph, we now suppose that there exists an $R$ such that $\D \Omega$ is a graph over $B(p,R)\cap T_{p}\D \Omega$ of height at most $\lambda_{0}$ for every $p\in \D \Omega$.  This follows from $C^{1}$ regularity of the boundary (and hence uniform continuity of the derivative. Let $\varphi_{p}:T_{p}\D \Omega\to \D \Omega$ be this function. We denote by $\nu(p)$ the inner normal to the boundary at $p$. We also consider $h_{p}=\<\varphi_{p}-p,\nu(p)\>$ the height of $\varphi(p)$ in the direction of $\nu(p)$.




Let $\Omega_{p,\lambda}$ denote the set $\{y\in \Omega, \<y,\nu(p)\>\geq \lambda\<p,\nu(p)\>\}$ for $\lambda\leq 1$ and $p\in \D \Omega$.

By $C^{1}$ regularity we have that for every $\alpha>0$ there is $\lambda_{0}$ which is independent of $p$, such that the cone $C_{\alpha,p}=\{y\mid \<p-y,\nu(p)\>\geq \alpha |y-p|\}$ is included in $\Omega_{p,\lambda_{0}}$.


\begin{figure}
\includegraphics[scale=0.4]{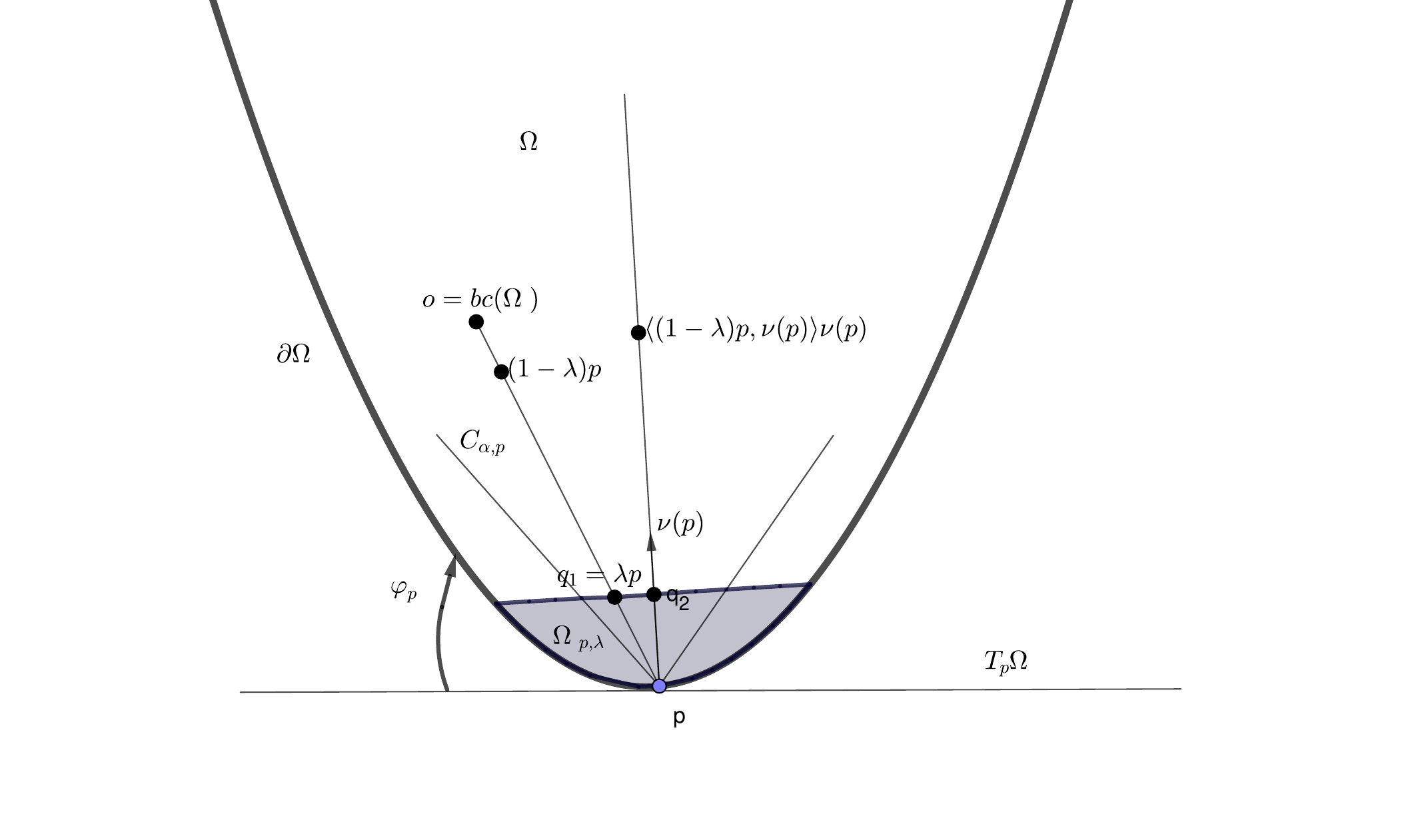}
\caption{\label{figure:aroundp} Notations of the sets involved for the local study around $p$}
\end{figure}

\begin{lem}
  \label{lem:upper_bound_1}
  Let $\lambda_{0}$ be small enough such that the cone $C_{\alpha,p}$ is contained in $\Omega_{p,\lambda_{0}}$ for some $\alpha\leq \<p/|p|,\nu(p)\>/2$, and every $p$. Then there is a number $C$ which depends on the $C^{1}$ norm of $\D \Omega$ such that for every $p$ and every $\lambda\leq \lambda_{0}$,
 \[\sigma(\lambda p)\leq C\sigma(p -\<(1-\lambda)p,\nu(p)\>\nu(p)).\]
\end{lem}
\proof
Let $q_{1}=\lambda p$ and let $q_{2}$ be in the normal direction to $T_p\partial \Omega$:
\[q_2=p-\<(1-\lambda)p,\nu(p)\>\nu(p).\]  

Let $c$ be a positive number such that for every $p\in \D \Omega$  $B(q_{1},\lambda c)\subset C_{\alpha,p}\cap \Omega_{p,\lambda}$ and $B(q_{2},(\lambda) c\subset C_{\alpha,p}\cap \Omega_{p,\lambda}$, and hence both balls are contained in $\Omega$. The existence of such a $c$ is given by the fact that we can fit a cone inside $\Omega_{p,\lambda_{0}}$.

Now for any $q\in \D \Omega\cap \Omega_{p,\lambda}$,
\begin{align*}
  |q-q_{1}|\leq |q-q_{2}|+|q_{1}-q_{2}|\\
  \intertext{And so}
  \frac{|q-q_{1}|}{|q-q_{2}|}\leq 1+\frac{|q_{1}-q_{2}|}{|q-q_{2}|},
\end{align*}
But because we can fit balls of radius $c(\lambda)$ around both $q_{i}$, and because $|q_{1}-q_{2}|\leq |\lambda| \sup_{p}|<\nu(p),p\> \leq |\lambda|$, it follows that there is a $C'$ such that
\[\frac{|q-q_{1}|}{|q-q_{2}|}\leq C'\;\;\mbox{and similarly}\;\;\frac{|q-q_{2}|}{|q-q_{1}|}\leq C'.\]

Now consider the half spheres $S^{+}_{\nu(p)}(q_i)\subset T_{q_{i}}\Omega$ ($i=1,2$) defined by
\[S^{+}_{\nu(p)}(q_i) = \set{v\in  T_{q_{i}}\; \norm{v}=1\;\mbox{and}\;\left\langle v,\nu(p) \right\rangle \leqslant 0}.\]
To compare the Busemann density at the points $q_1$ and $q_2$, we introduce the change of variables $S^{+}_{\nu(p)}(q_{1})\to S^{+}_{\nu(p)}(q_{2})$ obtained by first mapping $v$ to its projection (through $q_1$) on the boundary of $\Omega$ and then project it onto the Euclidean unit sphere centered at $q_{2}$. We now show that this change of variables is bi-Lipschitz.

Let $v_{1}$ and $v_{2}$ be two elements in $T_{q_{1}}$.  Let $\gamma:[0,\theta]\to S^{+}_{\nu(p)}(q_1)$ be the great arc connecting $v_{1}$ and $v_{2}$.  Let $\pi_{1}$ be the projection from the unit sphere centered at $q_{1}$, to $\D \Omega$ and let $\pi_{2}$ that from $q_{2}$. Consequently, if $\eta=\pi_{1}\circ \gamma$, we have
\[\gamma=\frac{\eta-q_{1}}{|\eta-q_{1}|}\]
and
\[\pi_{2}^{-1}\circ\pi_{1}\circ \gamma=\frac{\eta-q_{2}}{|\eta-q_{2}|}.\]
Then
\begin{align*}
  |\dot{\gamma}|&=\sqrt{\frac{|\dot{\eta}|^{2}}{|\eta-q_{1}|}+\frac{\<\dot{\eta},\eta-q_{1}\>^{2}|\eta-q_{1}|^{2}}{|\eta-q_{1}|^{6}}-2\frac{\<\eta-q_{2},\dot{\eta}\>^{2}}{|\eta-q_{1}|^{4}}}\\
  &=\frac{|\dot{\eta}|}{|\eta-q_{1}|^{2}}\sqrt{|\eta-q_{1}|^{2}-\<\dot{\eta}/|\dot{\eta}|,\eta-q_{1}\>^2}.
\end{align*}
It turns out that $\sqrt{|\eta-q_{1}|^{2}-\<\dot{\eta}/|\dot{\eta}|,\eta-q_{1}\>^2}$ is the norm of $\Pi_{\dot{\eta}^{\perp}}(q_{1}-\eta)$ ($\Pi_{\dot{\eta}^{\perp}}$ denotes the orthogonal projection on the hyperplane perpendicular to $\dot{\eta}^{\perp}$). The vector $q_{1}-\Pi_{\dot{\eta}^{\perp}}(q_{1}-\eta)$ is in the tangent hyperplane at $\eta$, hence $\Pi_{\dot{\eta}^{\perp}}(q_{1}-\eta)$ has a norm greater than $c\lambda$ (the size of the ball we can fit around $q_1$).
Now the key point is that $|(q_{1}-q_{2})|\leq |\lambda|$. 
We can apply a similar argument to bound $\frac{d}{dt}\pi_{2}^{-1}\circ \eta$, and
\begin{eqnarray*}
\frac{\abs{\frac{d}{dt}\pi_{2}^{-1}\circ \pi_{1}\circ\gamma}}{\abs{\dot{\gamma}}} & = & \frac{|\eta-q_{1}|^{2}}{|\eta-q_{2}|^{2}}\frac{|\Pi_{\dot{\eta}^{\perp}} (\eta-q_{2})|}{|\Pi_{\dot{\eta}^{\perp}}(\eta-q_{1})|}\\
 & \leqslant & \frac{|\eta-q_{1}|^{2}}{|\eta-q_{2}|^{2}}\frac{|\Pi_{\dot{\eta}^{\perp}} (\eta-q_{1})|}{|\Pi_{\dot{\eta}^{\perp}}(\eta-q_{1})|} + \frac{|\eta-q_{1}|^{2}}{|\eta - q_{2}|^{2}}\frac{|\Pi_{\dot{\eta}^{\perp}} (q_1 -q_{2})|}{|\Pi_{\dot{\eta}^{\perp}}(\eta-q_{1})|}\\
 & \leqslant & C'^2 + \frac{|\eta-q_{1}|^{2}}{|\eta - q_{2}|^{2}}\frac{| q_1-q_{2})|}{|\Pi_{\dot{\eta}^{\perp}}(\eta-q_{1})|}\\
 & \leqslant & C'^2+\frac{C'^2}{c}.
\end{eqnarray*}
We denote by $C$, the number $C'^2+\frac{C'^2}{c}$.
Hence
\[d(\pi_{2}^{-1}\circ\pi_{1}\circ \gamma(\theta),\pi_{2}^{-1}\circ \pi_{1}\circ \gamma(0))\leq C d(\gamma(\theta),\gamma(0)).\]
But by a similar argument
\[d(\pi_{2}^{-1}\circ\pi_{1}\circ \gamma(\theta),\pi_{2}^{-1}\circ \pi_{1}\circ \gamma(0))\geq C^{-1} d(\gamma(\theta),\gamma(0)).\]

Finally we open up $\sigma(q_{1})$:
\begin{align*}
& \sigma((1-\lambda) p)/\omega_{n}\\
& = \left(2\int_{S^{+}_{\nu(p)}(q_1)}\left(\frac{|\varphi_{q_{1}}(v)||\varphi_{q_{1}}(-v)|}{|\varphi_{q_{1}}(v)|+|\varphi_{q_{1}}(-v)|}\right)^{n}\:dv\right)^{-1}\\
& \leqslant  C^{2(n-1)}\left(\int_{S^{+}_{\nu(p)}(q_2)}\left(\frac{|\varphi_{q_{2}}(\pi_{2}^{-1}\circ \pi_{1}(v))||\varphi_{q_{2}}(\pi_{2}^{-1}\circ \pi_{1}(-v)|}{|\varphi_{q_{2}}(\pi_{2}^{-1}\circ \pi_{1}(v)|+|\varphi_{q_{2}}(\pi_{2}^{-1}\circ \pi_{1}(-v)|}\right)^{n}\:dv\right)^{-1}\\
& \leqslant  C^{3n-2}\left(\int_{S^{+}_{\nu(p)}(q_2)}\left( \frac{|\varphi_{q_{2}}(w)|\varphi_{q_{2}}(-w)|}{|\varphi_{q_{2}(w)}|+|\varphi_{q_{2}}(-w)|} \right)^{n}\:dw\right)^{-1}\\
& \leqslant   C^{3n-2}\sigma(q_{2})/\omega_{n}.
\end{align*}
\endproof

\begin{lem}
  \label{lem:sigma_cutoff_bound}Let $\Omega$ be a convex set, Let $x$ be a point on the boundary of $\Omega$, with inner normal $\nu(x)$.  Let $\varphi_{\lambda}$ be the projection of the unit ball centered at $x+\lambda \nu(x)$ to the boundary of $\Omega$. Assume that, for all $v\in S^{+}_{\nu(x)}(x)$,
$$d(x+\nu(x),\varphi_\lambda (v)) \leq C d(x+\nu(x),\varphi_\lambda(-v)).$$ Then 
  \[\sigma(x-\lambda \nu(x))\leq (C+1) |\Omega_{x,\lambda}|.\]
\end{lem}
\proof

We note that
\[|\Omega_{x,\lambda}|=\int_{S^{+}_{\nu(x)}}|\varphi_{\lambda}(v)|^{n}\:dv,\] 
\[\sigma(x-\lambda \nu(x))=\left(2\int_{S^{+}_{\nu(x)}(x)}\left(\frac{|\varphi(v)||\varphi(-v)|}{|\varphi(v)|+|\varphi(-v)|}\right)^{n}\:dv\right)^{-1}\]
and
\[\frac{|\varphi(v)||\varphi(-v)|}{|\varphi(v)|+|\varphi(-v)|}\geq \frac{|\varphi(v)||\varphi(-v)|}{(C+1)|\varphi(-v)|}\geq |\varphi(v)|/(C+1).\]
\endproof

\section{Centro-projective area for the Sobolev class}\label{sec:proofsobolev}

The point of this section is to prove the second main theorem:

\begin{thm}\label{thm:second}
  Let $\Omega$ be an open convex relatively compact subset of $\mathbb{R}^n$ such that the boundary has a regularity $W^{2,p}\cap C^{1} $ for some $p\geq n-1$. Then the Finsler volume growth entropy of $\Omega$ is maximal, equal to $n-1$
\end{thm}

Let us first remark that we have the Sobolev embedding theorem stating that the boundary has also regularity $\mathcal{C}^1$, so that lemma \ref{lem:upper_bound_1} can be applied in this situation.

This statement echoes the main result of \cite{bbvc11} for which the regularity assumption is $\mathcal{C}^{1,1}$ (hence $W^{2,\infty}$).  In fact the proof largely follows their proof: there are two crucial steps in the proof of the main theorem of \cite{bbvc11} that need to be worked out in this context. Precisely we want to show that
\begin{enumerate}
\item As a point $p\in \Omega$ approaches the boundary, the Busemann density (suitably renormalized) converges almost everywhere to $\frac{2a\sqrt{\det DG(p)}}{\left((1+a)\<\nu(p),p\>\right)^{\frac{n-1}{2}}}$. This will be achieved in lemma \ref{lem:ae_convergence}.
\item We then want to use this convergence to replace the volume of balls by the centro-projective area (i.e use the Lebesgue dominated convergence theorem). This requires a bound of the Busemann density by a dominating L$^1$ function. This will use the technology of maximal functions and will be achieved in lemma \ref{lem:integrable_bound}.
\end{enumerate}

Both of the lemmas \ref{lem:ae_convergence} and \ref{lem:integrable_bound} will follow from general considerations on the Busemann function (see paragraph \ref{section:busemann}).

\begin{lem}
  \label{lem:ae_convergence}
  For almost every $y\in B(x,R)$ 
\[\lim_{\lambda\to 1}\sigma(\lambda p)(1-\lambda)^{(n+1)/2}=\frac{\sqrt{k(p)}}{2^{\frac{n+1}{2}}\<p,\nu(p)\>}.\]
\end{lem}
\proof
We use Alexandrov's theorem on the almost everywhere second differentiability of convex functions \cite{Aleksandrov_twice_differentiable_convex} \cite{twice_differentiable_convex}, to yield for almost every $y\in B(x,R)$
\[\frac{1}{t^{2}}|h_{y}(t\theta)-t^{2}D^{2}h_{y}(0)(\theta,\theta)|\to 0 \]
as $t\to 0$, and $h_{y}:T_{y}\D \Omega\to \RR$ is the height of function so that $\{y+z-\nu(y)h(z)\mid z\in T_{y}\D \Omega\}$ is a subset of $\D \Omega$ (\emph{n.b.} $h_{y}(0)=0$).
The bilinear form $D^{2}h_{y}(0)$ is positive semidefinite, and given almost everywhere by the weak derivative \cite{nlrpotltheory}.  It has principal values $r_{1},\ldots r_{n-1}$ which are the principal curvatures, with principal directions $\tau_{1},\ldots, \tau_{n-1}$ which form an orthonormal basis.  For $\varepsilon\in \RR$ there is a $t_{0}$ such that for every $t<t_{0}$.
\[\frac{1}{t^{2}}|h_{y}(t\theta)-t^{2}D^{2}h_{y}(0)(\theta,\theta)|\leq \varepsilon, \]
so
\[[D^{2}h_{y}(0)-\varepsilon I](\theta,\theta)\leq h_{y}(t\theta)\leq [D^{2}h_{y}(0)+\varepsilon I](\theta,\theta).\]
This implies that the parabolas
\[P_{y,\varepsilon,\lambda}:=\{y+z-t\nu(y)\mid (D^{2}h_{y}(0)+\varepsilon I)(\zeta,\zeta)\leq t\leq \lambda,\]
satisfy
\[P_{y,\varepsilon,\lambda}\subset \Omega_{y,\lambda}\subset P_{y,-\varepsilon,\lambda}.\]
Finally we can proceed as in the proof of Proposition 2.8 in \cite{bbvc11}.
\endproof

We now consider a map $h_{x}:B(x,R)\subset T_{x}\D \Omega\to \D \Omega$.

Suppose that the restriction of $h_{x}$ to the line $t\mapsto y+t \theta$ is $W^{2,p}(\RR)$.  Then
\[ (y+t\theta)+h_{x}(y+t\theta)\nu(x)\]
\[=y+t\theta-\nu(x)h(y)-t \nu(x)Dh_{x}(y)\cdot \theta-\nu(x)\int_{0}^{t}\int_{0}^{\tau}\theta^{t}D^{2}h_{x}(y+s\theta)\theta\:ds\:d\tau.\]

\begin{lem}
  \label{lem:ae_density_points}
  Let $\{\theta_{i}\mid \in \mathbb{N}\}$ be a countable dense set of directions in $S^{n-2}\subset T_{x}\D \Omega$.  There is a set $E\subset B(x,R)$ of full measure such that for every $y\in E$ the map $f_{y,i}:t\mapsto h(y+t\theta_{i})$ is in $W^{2,p}(\RR)$, $f_{y,i}$ is twice differentiable at $t=0$ and
  \begin{equation}
    \label{eqn:second_derivative}
  \frac{d^{2}}{dt^{2}}|_{t=0}f_{y,i}=\theta_{i}^{t}D^{2}h(y)\theta_{i}.
\end{equation}
\end{lem}
\proof
This follows from the absolutely continuous on lines characterisation of Sobolev functions \emph{cf.} \cite{nlrpotltheory}.  For every direction we have a set $E'_{i}\subset B(x,R)\cap \theta_{i}^{\perp}$ such that the map $t\mapsto h(y+t\theta_{i}$ is in $W^{2,p}(\RR)$ by Fubini, and for almost every $t$ the second derivative is given by the weak second derivative of $h$ in the direction $\theta_{i}$.  Let $E_{i}$ be the set
\begin{multline*}
  E_{i}=\{y+t\theta_{i}\in B(x,R)\mid y\in E_{i}'\text{ and }s\mapsto h(y+s\theta_{i})\text{is twice differentiable}\\
\text{at $t$ with weak derivative }\theta_{i}^{t}D^{2}h(y+t\theta_{i})\theta_{i}\}.
\end{multline*}
Finally $E=\bigcap_{i\in \NN}E_{i}$.
\endproof

\begin{lem}
  \label{lem:integrable_bound}
  Assuming $\Omega$ is a domain of class $C^{1}\cap W^{2,p}$ for $p\geq (n-1)$, $p>1$ then there is a function $f\in L^{1}(\D \Omega)$ such that
  \[(1-\lambda)^{(n+1)/2}\sigma(\lambda p)\leq f(p)\]
  for all $\lambda\leq 1-\lambda_{0}$.
\end{lem}
\proof
Let $\theta_{1},\ldots \theta_{n-1}$ be an orthonormal basis for $T_{x}\Omega K$.  Define the set $E_{i}$ to be the set of $y$ such that $h_{x}|_{y+t\theta_{i}}$ is in $W^{2,p}(]a,b[)$.  Again $E_{i}$ is of full measure. Let $E=\bigcap_{i}E_{i}$.

  For every $i$ define the directional maximal function.
  \[M_{i}(f)(y)=\sup_{t\leq R/2}\frac{1}{2t}\int_{-t}^{t}|f(y+s\theta_{i})|\:ds\]
  This is bounded from $L^{p}(B(x,R))\to L^{p}(B(x,R/2))$ and for $1<p<\infty$, by Fubini's theorem and the boundedness of the usual maximal function on $L^{p}(\RR)\to L^{p}(\RR)$.

  For $n=2$ and $p=1$ this is bounded $L^{1}([x-R,x+R])\to L^{1,\infty}([x-R,x+R])$.

  We introduce several functions
  \[g_{i}=M_{i}(D^{2}h(\cdot)(\theta_{i},\theta_{i}))\]
  and $g=\max\{1,\max_{i} g_{i}\<\nu(x),\nu(y)\>\}$ all of which are in $L^{p}(B(x,R/2))$. For $y\in E$ such that $g(y)<\infty$
  \begin{align*}
    |h(y+t\theta_{i})-h(y)-tDh(y)\cdot \theta_{i}|&=\left|\int_{0}^{t}\int_{0}^{\tau}D^{2}h(y+s\theta_{i})(\theta_{i},\theta_{i})\:ds\:d\tau\right|\\
    &\leq 4t^{2}\frac{1}{4t^{2}}\int_{-t}^{t}\int_{-t}^{t}|D^{2}h(y+s\theta_{i})|\:ds\:d\tau\\
    &\leq 4t^{2}g_{i}(y)\leq 4t^{2}g(y).
  \end{align*}
  Let $\tau_{i}\in T_{y}\Omega$ be such that $\Pi_{T_{x}}(\tau_{i})=\theta_{i}$.  Let $\chi_{i,\lambda}$ be the intersection of the curve
  \[t\mapsto y+t\theta_{i}-t \nu(x) Dh(y)\cdot \theta_{i}+t^{2}g(y)\nu(x),\]
  and the plane
  \[L_{y,\lambda}=\{x\mid \<x,\nu(y)\>=\<y,\nu(y)\>-\lambda,\]
which we know to be in $\Omega$ for $t$ sufficiently small, and $L_{y,\lambda}$.  Then we know that 
\[|\eta_{i,\lambda}|=|\chi_{i,\lambda}+\nu(y)\lambda-\sqrt{\lambda}\tau_{i}/\sqrt{g(y)\<\nu(x),\nu(y)\>}|\leq \lambda \]
It follows that the convex hull of the points $\chi_{i,\lambda}$ has Hausdorff $n-1$ measure
given by
\[
  \sqrt{\frac{\lambda}{g(y)}}^{n-1} 4(\tau_{1}+\sqrt{g(y)/\lambda}\eta_{1,\lambda})\wedge\cdots \wedge (\tau_{n-1}+\sqrt{g(y)/\lambda}\eta_{n-1}).
\]
Assume $\lambda<1/(10 C(n)g(y))=\lambda_{0}(y)$, then this is greater than
\[\kappa\sqrt{\frac{\lambda}{g(y)}}^{n-1}(9/10)\]
where $\tau_{1}\wedge\cdots \wedge t_{n-1}=\kappa$.
Applying Lemma \ref{lem:sigma_cutoff_bound}
\begin{scriptsize}
\begin{align*}
\sigma(y&-\lambda \nu(y))\lambda^{(n+1)/2}\\
&\leq \frac{2\lambda^{(n+1)/2}}{\kappa} \left(\int_{0}^{\min\{\lambda,\lambda_{0}(y)\}}\sqrt{\frac{\lambda'}{g(y)}}^{n-1}\:d\lambda'+C(\Omega)\min\{\lambda-\lambda_{0}(y),0\}\frac{1}{c(n)g(y)^{n-1}}\right)^{-1}\\
&\leq C(n,\Omega,R)\sqrt{g(y)}^{n-1}
\end{align*}
\end{scriptsize}
if $\lambda<\lambda_{0}(y)$ or $\lambda_{0}(y)=\lambda_{0}$ and otherwise
\begin{align*}
&\leq C(n,\Omega,R)\frac{\lambda^{(n+1)/2}}{C(\Omega)\frac{\lambda_{0}(y)}{g(y)^{n-1}}+C(\Omega)(\lambda-\lambda_{0}(y))\frac{1}{g(y)^{n-1}C(n)10}}\\
&\leq C(n,\Omega,R) g(y)^{n-1}.
\end{align*}
Applying this inequality and Lemma \ref{lem:upper_bound_1} yields the result.

For $n=2$ we do not have the error terms, and get 
\[\sigma(y+\lambda \nu(y))\lambda^{(n+1)/2}\leq \sqrt{g(y)}.\]
But if $g\in L^{1,\infty}$ then $\sqrt{g}$ is integrable (on a bounded set).
\endproof

\begin{rmk}
It is worth commenting that although the centro-projective area would appear to require that $\Omega$ is $W^{2,p}$ for $p\geq (n-1)/2$, we can only get an integrable  bound in  $L^{1}(\D \Omega)$ for $p\geq n-1$.  This is to be expected as this is the natural exponent for a boundary of dimension $n-1$. But still the question arises whether this can be reduced.  It is possible that there is a higher integrability result for convex boundaries i.e. $p\geq n-1-\varepsilon$ implies $p\geq n-1$.
\end{rmk}

\proof[Proof of Theorem \ref{thm:second}] We can apply Lemmata \ref{lem:ae_convergence} and \ref{lem:integrable_bound} along with the dominated convergence theorem.  This allows us follow equation (26) in the proof of Theorem 3.1 in \cite{bbvc11}, and bring the limit into the integral, yielding the centro-projective area.

\bibliographystyle{alpha}
\bibliography{biblio}

\def\cprime{$'$} \def\dbar{\leavevmode\hbox to 0pt{\hskip.2ex \accent"16\hss}d}
  \def\cprime{$'$}
\begin{thebibliography}{{Tho}15}

\bibitem[Ale39]{Aleksandrov_twice_differentiable_convex}
A.~D. Alexandroff.
\newblock Almost everywhere existence of the second differential of a convex
  function and some properties of convex surfaces connected with it.
\newblock {\em Leningrad State Univ. Annals [Uchenye Zapiski] Math. Ser.},
  6:3--35, 1939.

\bibitem[{\'A}PT04]{hilbertvolume}
J.~C. {\'A}lvarez~Paiva and A.~C. Thompson.
\newblock Volumes on normed and {F}insler spaces.
\newblock In {\em A sampler of {R}iemann-{F}insler geometry}, volume~50 of {\em
  Math. Sci. Res. Inst. Publ.}, pages 1--48. Cambridge Univ. Press, Cambridge,
  2004.

\bibitem[BBV10]{bbvc11}
Gautier Berck, Andreas Bernig, and Constantin Vernicos.
\newblock Volume entropy of {H}ilbert geometries.
\newblock {\em Pacific J. Math.}, 245(2):201--225, 2010.

\bibitem[BCG95]{bcg2}
G.~Besson, G.~Courtois, and S.~Gallot.
\newblock Entropies et rigidit\'es des espaces localement sym\'etriques de
  courbure strictement n\'egative.
\newblock {\em Geom. Funct. Anal.}, 5(5):731--799, 1995.

\bibitem[BCP96]{twice_differentiable_convex}
Gabriele Bianchi, Andrea Colesanti, and Carlo Pucci.
\newblock On the second differentiability of convex surfaces.
\newblock {\em Geom. Dedicata}, 60(1):39--48, 1996.

\bibitem[Bus47]{busemannvolume}
Herbert Busemann.
\newblock Intrinsic area.
\newblock {\em Annals of Mathematics}, 48(2):234--267, 1947.

\bibitem[Cho94]{choucroun}
Francis Choucroun.
\newblock Arbres, espaces ultram\'etriques et bases de structure uniforme.
\newblock {\em Geom. Dedicata}, 53(1):69--74, 1994.

\bibitem[Cra09]{cramponentropydivisible}
Micka{\"e}l Crampon.
\newblock Entropies of strictly convex projective manifolds.
\newblock {\em J. Mod. Dyn.}, 3(4):511--547, 2009.

\bibitem[CV04]{entropyhilbertc2}
Bruno Colbois and Patrick Verovic.
\newblock Hilbert geometry for strictly convex domains.
\newblock {\em Geom. Dedicata}, 105:29--42, 2004.

\bibitem[CVV11]{entropypolytopes}
Bruno Colbois, Constantin Vernicos, and Patrick Verovic.
\newblock Hilbert geometry for convex polygonal domains.
\newblock {\em J. Geom.}, 100(1-2):37--64, 2011.

\bibitem[Eva98]{PDEs}
Lawrence~C. Evans.
\newblock {\em Partial differential equations}, volume~19 of {\em Graduate
  Studies in Mathematics}.
\newblock American Mathematical Society, Providence, RI, 1998.

\bibitem[HKM06]{nlrpotltheory}
Juha Heinonen, Tero Kilpel{\"a}inen, and Olli Martio.
\newblock {\em Nonlinear potential theory of degenerate elliptic equations}.
\newblock Dover Publications Inc., Mineola, NY, 2006.
\newblock Unabridged republication of the 1993 original.

\bibitem[IM01]{GFT}
Tadeusz Iwaniec and Gaven Martin.
\newblock {\em Geometric function theory and non-linear analysis}.
\newblock Oxford Mathematical Monographs. The Clarendon Press Oxford University
  Press, New York, 2001.

\bibitem[Kay67]{Kay}
David~C. Kay.
\newblock The ptolemaic inequality in {H}ilbert geometries.
\newblock {\em Pacific J. Math.}, 21:293--301, 1967.

\bibitem[Lew71]{Royden_Algebra_Space}
Lawrence~G. Lewis.
\newblock Quasiconformal mappings and {R}oyden algebras in space.
\newblock {\em Trans. Amer. Math. Soc.}, 158:481--492, 1971.

\bibitem[MS09]{saaranen}
Pertti Mattila and Pirjo Saaranen.
\newblock Ahlfors-{D}avid regular sets and bilipschitz maps.
\newblock {\em Ann. Acad. Sci. Fenn. Math.}, 34(2):487--502, 2009.

\bibitem[Pal68]{NLglobalanalysis}
Richard~S. Palais.
\newblock {\em Foundations of global non-linear analysis}.
\newblock W. A. Benjamin, Inc., New York-Amsterdam, 1968.

\bibitem[PT14]{handbookhilbert}
Athanase Papadopoulos and Marc Troyanov, editors.
\newblock {\em Handbook of {H}ilbert geometry}, volume~22 of {\em IRMA Lectures
  in Mathematics and Theoretical Physics}.
\newblock European Mathematical Society (EMS), Z\"urich, 2014.

\bibitem[{Tho}15]{tholozanhilbert}
N.~{Tholozan}.
\newblock {Entropy of Hilbert metrics and length spectrum of Hitchin
  representations in $\mathrm{PSL}(3,\mathbb{R})$}.
\newblock {\em ArXiv e-prints}, June 2015.

\bibitem[{Ver}12]{entropyapprox}
C.~{Vernicos}.
\newblock {Approximability of convex bodies and volume entropy in Hilbert
  geometry}.
\newblock {\em ArXiv e-prints}, July 2012.

\end{thebibliography}

\end{document}